\documentclass[a4paper,12pt]{amsart}
\usepackage{color}
\usepackage{amsmath,amssymb,amscd,amsfonts,amsthm,stmaryrd}
\usepackage[mathscr]{eucal}
\usepackage{xcolor}
\usepackage{graphicx}

\usepackage{nicefrac}

\usepackage{color}

\numberwithin{equation}{section}

\def\length{\operatorname{length}}

\def\diam{\operatorname{diam}}
\def\rank{\operatorname{rank}}
\def\pol{\operatorname{pol}}
\def\rad{\operatorname{rad}}
\def\id{\operatorname{id}}

\def\acts{\curvearrowright}
\def\D{\partial}

\def\R{\mathbb R}
\def\Z{\mathbb Z}

\def\N{\mathbb N}
\def\Ha{\mathcal H}

\newcommand{\F}{F}
\newcommand{\G}{\Theta}
\newcommand{\Fi}{\F_\infty}
\newcommand{\Gi}{\G_\infty}

\def\al{\alpha}

\def\ka{\kappa}

\def\eps{\epsilon}
\def\ga{\gamma}
\def\Ga{\Gamma}

\def\la{\lambda}
\def\La{\Lambda}

\def\om{\omega}

\def\si{\sigma}
\def\Si{\Sigma}

\def\tits{\partial_{T}}
\def\geo{\partial_{\infty}}
\def\etits{\partial^{ess}_{T}}

\newcommand{\bI}{{\mathbf I}}
\newcommand{\bZ}{{\mathbf Z}}
\newcommand{\M}{{\mathbf M}}

\newcommand{\cD}{{\mathscr D}}
\newcommand{\cF}{{\mathscr F}}

\newcommand{\cP}{{\mathscr P}}

\newcommand{\on}{\:\mbox{\rule{0.1ex}{1.2ex}\rule{1.1ex}{0.1ex}}\:}

\def\acts{\curvearrowright}

\def\Isom{\mathop{\hbox{Isom}}}
\def\Stab{\mathop{\hbox{Stab}}}
\def\spt{\operatorname{spt}}

\def\<{\langle}
\def\>{\rangle}

\newcommand{\bb}[1]{\llbracket #1\rrbracket} 

\newcommand{\loc}{\text{\rm loc}}

\theoremstyle{plain}
\newtheorem{thm}{Theorem}[section]
\newtheorem{lem}[thm]{Lemma}
\newtheorem{prop}[thm]{Proposition}
\newtheorem{cor}[thm]{Corollary}
\newtheorem{slem}[thm]{Sublemma}

\newtheorem{introthm}{Theorem}

\newtheorem{introcor}{Corollary}

\newtheorem*{conj}{Conjecture}
\newtheorem*{quest}{Question}

\newtheoremstyle{named}{}{}{\itshape}{}{\bfseries}{.}{.5em}{\thmnote{#3} #1}
\theoremstyle{named}
\newtheorem*{namedlemma}{Lemma}

\theoremstyle{definition}
\newtheorem{dfn}[thm]{Definition}

\theoremstyle{remark}

\newcommand{\bcl}{\begin{claim}}
\newcommand{\ecl}{\end{claim}}
\newcommand{\bcor}{\begin{cor}}
\newcommand{\ecor}{\end{cor}}
\newcommand{\bdfn}{\begin{dfn}}
\newcommand{\edfn}{\end{dfn}}
\newcommand{\ben}{\begin{enumerate}}
\newcommand{\bit}{\begin{itemize}}
\newcommand{\blem}{\begin{lem}}
\newcommand{\bslem}{\begin{slem}}
\newcommand{\bprop}{\begin{prop}}
\newcommand{\bthm}{\begin{thm}}
\newcommand{\een}{\end{enumerate}}
\newcommand{\eit}{\end{itemize}}
\newcommand{\elem}{\end{lem}}
\newcommand{\eslem}{\end{slem}}
\newcommand{\eprop}{\end{prop}}
\newcommand{\ethm}{\end{thm}}

\usepackage{hyperref}



\begin{document}

\title[Higher rank]{CAT(0) spaces of higher rank}
\author{Stephan Stadler}

\newcommand{\Addresses}{{\bigskip\footnotesize
\noindent Stephan Stadler,
\par\nopagebreak\noindent\textsc{Max Planck Institute for Mathematics, Vivatsgasse 7, 53111 Bonn, Germany}
\par\nopagebreak
\noindent\textit{Email}: \texttt{stadler@mpim-bonn.mpg.de}

}}



\begin{abstract}
Ballmann's Rank Rigidity Conjecture predicts that a CAT(0) space of higher rank with a geometric group action
is {\em rigid} --
isometric to a Riemannian symmetric space, a Euclidean building, or splits as a direct product. 
We confirm this conjecture in rank 2.
For CAT(0) spaces of higher rank $n\geq 2$ we prove rigidity if the space contains a periodic $n$-flat and every complete geodesic lies in some $n$-flat.
This does not require a geometric group action.
\end{abstract}

\maketitle

\tableofcontents

\section{Introduction}
\subsection{Main results}
A {\em Hadamard manifold} is a simply connected Riemannian manifold of non-positive sectional curvature.
Its {\em rank} is the maximal dimension of  an isometrically embedded Euclidean space, a so-called {\em flat}.
The structure of manifolds of non-positive curvature and higher rank was clarified in the eighties, in the series of papers \cite{BBE_structure,BBS_structure},
culminating in Ballmann's celebrated Rank Rigidity Theorem \cite{B_higher}: The universal covering of a non-positively curved manifold of finite volume 
and higher rank is either a Riemannian symmetric space or splits as a direct product. 

The modern theory of non-positive curvature originates from Gromov's seminal paper \cite{G_hyp}.
Nowadays,
CAT(0) spaces -- the synthetic versions of Hadamard manifolds -- play an important role in mathematics beyond geometry and topology.
They appear in many different fields such as geometric group theory,
representation theory, arithmetic and optimization.
For an account of the huge amount of literature on metric spaces with synthetic curvature bounds, see the bibliographies in \cite{AKP, Ballmann, BH}. 
The following conjecture, formulated by Ballmann and Buyalo in \cite{BaBu_periodic}, is the main motivation for the present paper.

\begin{conj}[Rank Rigidity]
Let $X$ be a locally compact geodesically complete CAT(0) space. Suppose that $\Ga$ is a group of isometries of $X$
with limit set $\La(\Ga)=\geo X$. Then $\diam(\tits X)=\pi$ implies that $X$ is a Riemannian symmetric space or a Euclidean
building of rank at least 2, or that $X$ non-trivially splits as a direct product.
\end{conj}
The condition on the diameter is equivalent to saying that any complete geodesic bounds a flat half-plane.
If a group $\Ga$ acts on a CAT(0) space $X$, then its {\em limit set} $\La(\Ga)$ is the set of accumulation points  
of a $\Ga$-orbit in the ideal boundary $\geo X$. The main example of group actions on CAT(0) spaces with full limit set are {\em geometric actions}, meaning
cocompact and properly discontinuous actions by isometries. In the smooth case, namely for Hadamard manifolds, it is enough to assume finite covolume instead of 
cocompactness. 

The conjecture is known for Hadamard manifolds  if
$\Ga$ acts properly discontinously and with finite covolume \cite{B_higher,BS_higher}, or if $\Ga$ satisfies the duality condition \cite{ballmannbook, EH_diff}.
It is also known for homogeneous Hadamard manifolds \cite{Heber}.
For singular spaces, the conjecture is known to hold

\begin{itemize} 
	\item for a 2-dimensional simplicial complex with a piecewise
smooth CAT(0) metric and a geometric group action \cite{BaBr_orbi};
\item for a 3-dimensional simplicial complex with a piecewise Euclidean CAT(0) metric and a geometric group action \cite{BB_rr};
\item for a finite dimensional CAT(0) cube complex with a geometric group action 
\cite{CS_rr}\footnote{For cube complexes, the statement of rank rigidity takes the form of a splitting result, since
neither symmetric spaces nor Euclidean buildings carry a piecewise cubical CAT(0) metric \cite{Leeb}.};
\item for a locally compact and geodesically complete CAT(0) space such that its isometry group does not fix a point at infinity and
which supports a cocompact isometric action by an amenable locally compact group \cite{CM_bieber}. 
\item for a space of rank 2 which supports a geometric group action by a group $\Ga$ and whose
Tits boundary contains a proper closed $\Ga$-invariant subset \cite{Ricks_1D}.
\end{itemize}

Note that the interesting case in the next to last item concerns non-discrete groups since the only locally compact geodesically complete CAT(0)
space which admits a cocompact isometric action by an amenable group is the flat Euclidean space \cite[Corollary~C]{AB_amenable}.
\medskip  

We prove:

\begin{introthm}\label{thm_mainA}
Let $X$ be a locally compact and geodesically complete CAT(0) space of rank $n\geq 2$. 
Suppose that every complete geodesic in $X$ lies in an $n$-flat. 
If $X$ contains a periodic $n$-flat, then $X$ is a Riemannian symmetric space or a Euclidean building, or  $X$ splits as a direct product.
\end{introthm}

\begin{introcor}\label{cor_mainA}
Let $X$ be a locally compact and geodesically complete CAT(0) space of rank $n\geq 2$. 
Suppose that every complete geodesic in $X$ lies in an $n$-flat. 
If $X$ contains a periodic $n$-flat, then $X$ has a unique decomposition into irreducible factors
\[X=\R^j\times X_1\times\ldots \times X_k\times Y_1\times\ldots \times Y_l\times M_1\times\ldots\times M_m\]
where $X_a$ are rank 1 spaces, $Y_b$ are irreducible Euclidean buildings and $M_c$ are
irreducible Riemannian symmetric space of higher rank.
\end{introcor}

Note that in the smooth case, rank $n$ is expressed by the property that every complete geodesic admits $n$ linearly independent parallel Jacobi fields.
In the presence of a geometric group action these Jacobi fields can then be ``integrated'' to  $n$-flats \cite[Section~IV.4]{ballmannbook}.
Hence, in the context of metric spaces, the assumption on the existence of $n$-flats in Theorem~\ref{thm_mainA} is natural.
Recall that a flat  is {\em periodic}, if its stabilizer contains a subgroup acting geometrically. 
We point out that Theorem~\ref{thm_mainA} does not assume any symmetries besides the existence of a single periodic $n$-flat.
Also, if $X$ is a rank $n$ symmetric space or  Euclidean building, and $\Ga$ acts cocompactly on $X$, then $X$
contains a periodic $n$-flat \cite[Theorem~8.9]{BaBr_orbi}, \cite[Lemma~8.3]{Mos_rigidity}, \cite[Theorem~2.8]{PR_cartan}.
We emphasize that some symmetry is required, since a geodesically complete CAT(0) space of higher rank and whose  Tits boundary has diameter $\pi$
does not have to be a product or a Euclidean building \cite{BB_rr}.
Admittedly, it is unknown if a locally compact and geodesically complete CAT(0) space of rank $n$, where every complete geodesic lies in an $n$-flat,
and which supports a cocompact group action, has to contain a periodic $n$-flat. 

For spaces of rank 2, we prove:

\begin{introthm}\label{thm_mainB}
Let $X$ be a locally compact and geodesically complete CAT(0) space of rank $2$. 
Suppose that the diameter of its Tits boundary $\tits X$ is equal to $\pi$. 
If $X$ admits a geometric group action $\Ga\acts X$, then $X$ is a Riemannian symmetric space or a Euclidean building, or  $X$ non-trivially splits as a direct product.
\end{introthm}

This settles the Rank Rigidity Conjecture in rank 2, at least for spaces with a geometric group action.

We point out that we do not deduce Theorem~\ref{thm_mainB} from Theorem~\ref{thm_mainA} as we were unable to a priori show the existence of
a periodic 2-flat under the assumptions of Theorem~\ref{thm_mainB}. 
While the overall strategies for Theorems~\ref{thm_mainA} and~\ref{thm_mainB} are similar, the possible lack of a periodic 2-flat
in  Theorem~\ref{thm_mainB} leads to technical difficulties and the need for geometric measure theory methods.

\subsection{Strategy and predecessors}

Our proofs rely on the study of the asymptotic geometry of CAT(0) spaces, in particular the Tits boundary as introduced by Gromov in
\cite{BGS_npc}. This approach is inspired by Mostow's Strong Rigidity \cite{Mos_rigidity}: Two irreducible locally symmetric spaces of noncompact type
and rank $\geq 2$ with isomorphic fundamental groups are isometric up to scaling. Gromov generalized Mostow's theorem to the case where only one of the two spaces is assumed to be locally symmetric \cite{BGS_npc}, see also \cite{BE_fund,Eb_rig}. While in Mostow's original proof, the Tits metric was only implicitly present, Gromov's argument relies on a precise study of Tits boundaries. Putting complete focus on the asymptotics of a non-positively curved manifold, Burns and Spatzier introduced the notion of topological Tits buildings \cite{BS_buildings} and were able to give an alternative proof of Ballmann's
Rank Rigidity Theorem \cite{BS_higher}, adapting  Gromov's arguments. The counterpart of Mostow's theorem for Euclidean buildings was proved by Prasad \cite{Pr_lattices}.
It is natural to wonder if asymptotic rigidity extends even further, from Hadamard manifolds to CAT(0) spaces, quasi from Riemannian geometry to Alexandrov geometry.
This is indeed true and the content of Leeb's Rigidity Theorem.

\bthm[{\cite[Main Theorem]{Leeb}}]\label{thm_leeb}
Let $X$ be a locally compact, geodesically complete CAT(0) space. If $\tits X$ is a connected thick irreducible spherical building,
then $X$ is either a symmetric space or a Euclidean building.
\ethm

The proofs of Theorem~\ref{thm_mainA} and \ref{thm_mainB} use Leeb's result.
In order to show that the Tits boundary is a spherical building, we rely on Lytchak's Rigidity Theorem.

\bthm[{\cite[Main Theorem]{Ly_rigidity}}]\label{thm_lytchak}
Let $Z$ be a finite dimensional geodesically complete
CAT(1) space. If $Z$ has a proper closed subset $A$ containing with each
point all of its antipodes, then $Z$ is
a spherical join or a spherical building.
\ethm

The strategy to approach rank rigidity in this way is well-known.

There are two main issues preventing one to implement this strategy.
First, we have to ensure geodesic completeness of the Tits boundary $\tits X$.
And secondly, we need to find
a proper closed set $A\subset\tits X$ which is also closed under taking antipodes.

We overcome these obstacles as follows.

Under the assumptions of Theorem~\ref{thm_mainA} any point at infinity is contained in a round $(n-1)$-sphere.
This implies  geodesic completeness by a topological argument using local homology \cite{Kleiner, LS_affine}. 
For Theorem~\ref{thm_mainB} we first work with a geodesically complete subspace $\etits X$ called the {\em essential Tits boundary}
which consists of the union of all simple closed geodesics in $\tits X$. We then show that $\etits X$ is a
spherical building and, in a second step, that it has to agree with the ordinary Tits boundary.

The second issue is more significant. We prove that the subset $O\subset \tits X$
of {\em regular points} -- points which have a neighborhood isometric to an open set in a round $(n-1)$-sphere --
is non-empty and closed under taking antipodes. As a consequence, either $\tits X$ is equal to a round $(n-1)$-sphere,
or $A=\tits X\setminus O$ is subject to Lytchak's Theorem~\ref{thm_lytchak}.

\subsection{How to find regular points}
The existence of regular points is the key to both theorems. In this section we provide an informal discussion of the arguments that yield regular points in Theorems~\ref{thm_mainA} and \ref{thm_mainB}.
\subsubsection{Regular points at infinity of 2-dimensional spaces}
Let us try to explain the idea for a locally compact CAT(0) space $X$ of dimension 2. If the Tits boundary of our space has finite diameter, but
does not contain a regular point, then it basically branches everywhere. Intuitively this should contradict the local compactness of our space.
Now suppose that $F\subset X$ is a periodic 2-flat. If its Tits boundary $\tits F$ does not contain a regular point, then 
$\tits F$ contains a dense set of branch points. Using periodicity, it is not hard to see that 
each branch point at infinity yields a flat half-plane $H$ orthogonal to $F$. Since $X$ has dimension 2, any two such half-planes have to be orthogonal to each other.
By periodicity of $F$, there is a compact set $K\subset F$ and an infinite sequence of pairwise orthogonal flat half-planes which intersect $K$. This does indeed contradict local compactness.
\subsubsection{The rank 2 case of Theorem~\ref{thm_mainA}}
The actual proof of Theorem~\ref{thm_mainA} in rank 2 is not much harder then the 2-dimensional case sketched above. 
Again, from periodicity, one obtains families of orthogonal flat half-planes. 
The key point is that these are again pairwise orthogonal (Lemma~\ref{lem_proj}). This orthogonality can be explained as follows.
Suppose $F$ is a 2-flat in our rank 2 space $X$. For simplicity we assume that there is no 2-flat parallel to $F$ in $X$. 
Let $H^\pm$ be two flat half-planes whose boundaries $\D H^\pm\subset F$ are not parallel.
Let $\eta^\pm\in\tits X$ denote the centers of the half-circles $\tits H^\pm\subset\tits X$.
We think of $\D H^\pm$ to be the nearest point projections of $\eta^\pm$ to $F$.
This is justified by the orthogonality of $H^\pm$ to $F$ which ensures that $\D H^\pm$ minimizes the relative distance to $\eta^\pm$ on $F$. 
 Now if the angle between $H^-$ and $H^+$
would be strictly less than $\frac{\pi}{2}$, then by a first variation argument we could strictly decrease this relative distance (Busemann level) 
and arrive at a contradiction.
(See Lemma~\ref{lem_illu} for more details.)  

The  case of rank larger than 2 is naturally  a little more complicated but follows the same scheme.

\subsubsection{Regular points in Theorem~\ref{thm_mainB}}
The setting of Theorem~\ref{thm_mainB} awaits with a new difficulty. Since we do not assume the existence of a periodic flat, 
we do not get orthogonal flat half-planes right away. However, a geometric group action still has a periodic 1-flat, an axis $g$.
Let us assume that $g$ lies in a 2-flat $F$. 
From the assumption on the Tits diameter, every branch point $\xi$ in $\tits F$ together with its antipode $\hat \xi$ yields a sequence of Tits geodesics
in $\tits X\setminus \tits F$ joining $\xi$ to $\hat \xi$ and whose lengths converge to $\pi$.
If we are able to extract a limit geodesic which does not run in $\tits F$, then we can produce an orthogonal  flat half-plane corresponding to the branch point
$\xi$ and argue as before.
While everything up to now can be done purely geometrically, this part needs some analysis. 
Intuitively, a geodesic $c$ in $\tits X$ of length approximately $\pi$ which  joins two antipodal points in $\tits F$
should bound a surface in $X$ which resembles a flat half-plane orthogonal to $F$.
We find such a surface by means of geometric measure theory, namely by solving an asymptotic Plateau problem relative to $F$.
More precisely, we find a minimal surface $S$ whose boundary lies in $F$ and which is asymptotic to $c$.
The relative minimality ensures that locally near every point in $\D S$ the surface has a definite distance from $F$.
This will serve as a substitute for orthogonality. Since the length of $c$ is close to $\pi$, the surface $S$ has almost Euclidean area growth.
We also show that the boundary $\D S$ has to intersect the axis $g$. Using the axial isometries we can then produce the required flat half-plane as 
a limit of relative minimizers, appropriately translated by the axial isometry.

\subsection{Open questions}

Recall (a special case of) Gromov's Flat Closing Problem:

\begin{quest}[{Flat Closing \cite[6.B{\tiny 3}]{G_asym}}]
Let $X$ be a locally compact geodesically complete CAT(0) space with a geometric group action $\Ga\acts X$. 
Suppose that $X$ contains an $n$-flat. Does this imply that $X$ contains a $\Ga$-periodic $n$-flat?
\end{quest}
In general, the answer is expected to be no,  cf.~\cite{KP_aperiodic, SW_periodic}. Yet, the Flat Closing Problem has a positive answer in a couple of cases, see the introduction to \cite{KP_aperiodic}.
\medskip

As mentioned above, Flat Closing holds for Riemannian symmetric spaces and Euclidean buildings. Theorem~\ref{thm_mainA} shows that
Flat Closing for rank $n$ spaces where every complete geodesic lies in an $n$-flat is equivalent to  Rank Rigidity  (in the cocompact case).

Since Theorem~\ref{thm_mainB} proves Rank Rigidity in rank 2 it also confirms Flat Closing for rank 2 spaces of Tits-diameter $\pi$.
We point out that if a group $\Ga$ acts geometrically on a CAT(0) space $X$, then $\Ga$ contains a hyperbolic element, 
namely one that acts as an axial isometry on $X$ \cite[Theorem~11]{Sw_cut}. Formally, this corresponds to the $n=1$ case of Flat Closing.
Our proof of Theorem~\ref{thm_mainB}
relies on this fact.

In view of Theorem~\ref{thm_mainA}, we would like to state the following interesting special case of Gromov's question.

\begin{conj}[{Higher Rank Flat Closing}]
Let $X$ be a locally compact geodesically complete CAT(0) space whose Tits-diameter is $\pi$ and which supports a geometric group action $\Ga\acts X$. 
If $X$ has rank $n$, then $X$ contains a periodic $n$-flat.
\end{conj}

By the existence of asymmetric Euclidean buildings, we genuinely expect Theorem~\ref{thm_mainA} to treat a more general setting than Rank Rigidity.
Actually, we obtain Theorem~\ref{thm_mainA} as  a consequence of the following.

\begin{introthm}[Theorem~\ref{thm_open}]\label{thm_mainC}
Let $X$ be a locally compact CAT(0) space of rank $n\geq 2$.  If $X$ contains a periodic $n$-flat $F$, then $\si=\tits F$ contains a dense subset
$U\subset \si$ which is open in $\tits X$.
\end{introthm}

This theorem aims at to following general conjecture:

\begin{conj}[Lytchak]
Let $X$ be a locally compact CAT(0) space with a geometric group action. 
Then $\tits X$ contains a relatively compact open subset.
\end{conj}

Again, if Gromov's Flat Closing would hold for higher rank CAT(0) spaces, then Theorem~\ref{thm_mainC} would confirm Lytchak's conjecture.

\subsection{Organization}

In Section~\ref{sec_pre}, we agree on notation and recall background on the geometry of spaces with upper curvature bounds.
In Section~\ref{sec_buildings}, we recall metric characterizations of spherical and Euclidean buildings. 
The central content of this section is a recognition theorem for buildings (Theorem~\ref{thm_rec}) which is later used in the proof of Theorem~\ref{thm_mainA}.
The short Section~\ref{sec_illu} is included for the convenience of the reader. It illustrates the basic idea behind discreteness for branch points at infinity in 
a simplified setting. 
We begin Section~\ref{sec_rank_n} by recalling how branch points at infinity of periodic flats lead to flat half-planes orthogonal to the periodic flat.
The goal of this section is then to  show that a periodic $n$-flat in a space of rank $n$ leads to a non-empty set of regular points at infinity
(Theorem~\ref{thm_open}). An important ingredient is Proposition~\ref{prop_halfspace} which constructs a flat 3-dimensional half-space from a sequence of flat 
half-planes. We conclude this section with the proof of Theorem~\ref{thm_mainA}.
Section~\ref{sec_rank_2} builds up to a proof of Theorem~\ref{thm_mainB} assuming our key technical result, the \hyperref[lem_key_tech]{Half-Plane Lemma}.
Appenix~\ref{sec_app} is devoted to geometric measure theory with the aim to provide a proof of the \hyperref[lem_key_tech]{Half-Plane Lemma}.
We also prove some folklore results, including monotonicity and volume rigidity of minimizing currents in CAT(0) spaces.
However, our approach is to find a streamlined path to  the \hyperref[lem_key_tech]{Half-Plane Lemma} and we do not try to prove
required ingredients in their most general form.

\subsection{Acknowledgments}
It's my pleasure to thank Alexander Lytchak for carefully reading a late version of this paper,
pointing out several mistakes,
and making insightful comments that led to substantial improvements. I also want to thank Anton Petrunin for a helpful discussion on
Proposition~\ref{prop_mon}. 
I was supported by DFG grant SPP 2026.

\section{Preliminaries}\label{sec_pre}

General references for this section are \cite{AKP, Ballmann, BH, KleinerLeeb}.

\subsection{Metric spaces}
Euclidean $n$-space with its flat metric will be denoted by $\R^n$. The unit sphere $S^{n-1}\subset\R^n$ equipped with the induced metric will be referred to as a
{\em round sphere}. Its intersection with a half-space $\R^{n-1}\times[0,\infty)$ is called a {\em round hemisphere}.
We denote  the distance between two points $x$ and $y$ in a metric space $X$ by $|x,y|$.
If $A\subset X$ denotes a subset, then $|x,A|$ refers to the greatest lower bound for distances from points in $A$ to $x$.
For subsets $A, A'\subset X$ we denote the Hausdorff (pseudo-) distance by $|A,A'|_H$.  
For $x\in X$ and $r>0$, we denote by $B_r(x)$ and $\bar B_r(x)$ the open and closed $r$-ball around $x$, respectively.
Similarly, $N_r(A)$ and $\bar N_r(A)$ denote the open and closed $r$-neighborhood of a subset $A\subset X$, respectively.
Moreover, $S_r(x)$ denotes the $r$-sphere around $x$ and by $\dot B_r(x)$ we denote the punctured $r$-ball $B_r(x)\setminus\{x\}$.
A \emph{geodesic}
is an isometric embedding of an interval. It is called a {\em geodesic segment}, if it is compact.
The {\em endpoints} of a geodesic segment $c$ are denoted by $\D c$.
A geodesic segment $c$ {\em branches} at an endpoint $y\in \D c=\{x,y\}$, if there are geodesics $c^\pm$ starting in $x$ which  strictly contain $c$
and such that $c^-\cap c^+=c$.

 A \emph{triangle} is a union of three geodesics connecting three points.
$X$  is \emph{a  geodesic metric space} if
any pair of   points of $X$
is connected by a geodesic.
It is \emph{geodesically complete} if every geodesic segment is contained in a complete geodesic.

\subsection{Spaces with an upper curvature bound}

For $\kappa \in \R $, let $D_{\kappa}  \in (0,\infty] $ be the diameter of the  complete, simply connected  surface $M^2_{\kappa}$
of constant curvature $\kappa$. A complete  metric space is called a \emph{CAT($\kappa$) space}
if any pair of its points with distance less than $D_{\kappa}$ is connected by a geodesic and if
 all triangles
with perimeter less than $2D_{\kappa}$
are not thicker than
 the \emph{comparison triangle} in $M_{\kappa} ^2$. In particular, geodesics between points of distance less than $D_{\kappa}$
are unique. Hence, $X$ is a CAT($\ka$) space, then we can define for every $p\in X$ and every subset $A\subset B_{D_{\ka}}(p)$ 
the {\em geodesic cone $C_p(A)$}.

For any CAT($\kappa$) space  $X$,
the angle between each pair of geodesics starting at the same point
is well defined. 
 Let $x,y,z$ be three points at pairwise distance less than $D_{\kappa}$ in a CAT($\kappa$) space $X$.
Whenever $x\neq y$, the geodesic between $x$ and $y$ is unique and will be denoted
by $xy$.   For $y,z \neq x$, the angle at $x$ between $xy$ and $xz$
will be denoted by $\angle_x(y,z)$.

The \emph{space of directions}
 $\Sigma _x X$ at  $x\in X$,  equipped with the angle metric,
is a CAT(1) space.
The Euclidean cone
over $\Sigma _x$ is a CAT(0) space.
It is denoted by $T_x X$ and called
the \emph{tangent space} at $x$ of $X$. Its tip will be denoted by $o_x$.
There are natural {\em logarithm maps}
\[\log_x:\dot B_x(D_\ka)\to\Si_x X,\quad\quad\log_x:X\to T_x X.\]

 A natural notion of dimension $\dim (X)$ for a CAT($\ka$) space $X$ was
  introduced by
Kleiner in \cite{Kleiner}.
It is equal to the supremum of topological dimensions of compact subsets of  $X$, satisfies 
\[\dim (X)= \sup _{x\in X} \{ \dim (\Sigma _xX) +1 \}\] 
and coincides with the homological dimension of $X$.

\subsection{CAT(1) spaces}

For two CAT(1) spaces $Z_1$ and $Z_2$ we denote by $Z_1\circ Z_2$ their {\em spherical join}.
It is a CAT(1) space of diameter $\pi$.
Note that every round subsphere $\si'$ in a round sphere $\si$ yields a join decomposition $\si=\si'\circ\si''$
where $\dim(\si)=\dim(\si')+\dim(\si'')+1$. Also, every round hemisphere $\tau$ decomposes as $\tau=\D\tau\circ\zeta$
where $\zeta$ denotes the center of $\tau$. A CAT(1) space $Z$ is called {\em irreducible}, if it does not admit a non-trivial
spherical join decomposition.

 A subset $C$ in a CAT(1) space $Z$ is called {\em $\pi$-convex}
if for any pair of points $x, y\in Z$ at distance less than $\pi$ the unique geodesic $xy$ is contained in $C$.
If $C$ is closed, then it is CAT(1) with respect to the induced metric.
For instance, a ball of radius at most $\frac{\pi}{2}$ is $\pi$-convex.
Let $C\subset Z$ be a closed convex subset with radius $\rad(C)\geq \pi$. Then we define the set of {\em poles} for $C$ by
\[\pol(C):=\{\eta\in Z|\ d(\eta,\cdot)|_{C}\equiv\frac{\pi}{2}\}.\]
If $\diam(C)>\pi$, then $C$ has no pole. If $\diam(C)=\rad(C)=\pi$, then $\pol(C)$ is closed and convex.
The convex hull of $C$ and $\pol(C)$ is canonically isometric to $C\circ\pol(C)$.
In case $C$ consists of a pair of {\em antipodes} $\xi^\pm$, $d(\xi^-,\xi^+)=\pi$, then the convex hull of
$\{\xi^-,\xi^+\}$ and $\pol(\{\xi^-,\xi^+\})$ is isometric to a spherical suspension of $\pol(\{\xi^-,\xi^+\})$.

A subset  $A\subset Z$ is called {\em spherical}, if it 
embeds isometrically into a round sphere. 
A {\em (spherical) $n$-lune of angle $\al\in[0,\pi]$} in a CAT(1) space is a closed convex subset $L$ isometric to 
$S^{n-1}\circ[0,\al]$.
Bigons in CAT(1) spaces lead to spherical lunes \cite[Lemma~2.5]{BB_diam}.
We will use the following more general result which follows immediately from \cite[Lemma~2.5]{BB_diam} and \cite[Lemma~4.1]{Ly_rigidity}.

\blem[Lune Lemma]\label{lem_lune}
Suppose that $\tau^\pm$ are round $n$-hemispheres in a CAT(1) space $Z$.
If $\D \tau^-=\D \tau^+$, then $\tau^-\cup \tau^+$ is a spherical subset in $Z$. 
More precisely, if $\angle(\tau^-,\tau^+)\geq \pi$, then $\tau^-\cup \tau^+$ is an $n$-sphere in $Z$;
and if  $\angle(\tau^-,\tau^+)=\al< \pi$, then $\tau^-\cup \tau^+$ bounds a spherical lune of angle $\al$.
\elem

\bdfn\label{def_reg_point}
Let $Z$ be a CAT(1) space of dimension $n-1$. We call a point $\xi\in Z$ {\em regular}, if
there exists  a radius $r>0$ such that $B_r(\xi)\subset Z$ is isometric to an $r$-ball in $S^{n-1}$.
\edfn

The following criterion for geodesic completeness is well-known \cite[Lemma~2.1]{BL_building}.

\blem\label{lem_gc}
	Let $Z$ be a CAT(1) space of dimension $n$. If every point $\xi\in Z$ is contained in a round $n$-sphere in $Z$,
	then $Z$ is geodesically complete.
\elem

\subsection{CAT(0) spaces}

The {\em ideal boundary} of a CAT(0) space $X$, equipped with the cone topology, is denoted by $\geo X$. 
If $X$ is locally compact, then $\geo X$ is compact. If $X$ is a Hadamard $n$-manifold, then $\geo X$ is homeomorphic to an $(n-1)$-sphere.  
The {\em Tits boundary} of $X$ is denoted by $\tits X$, it is the ideal boundary equipped with 
the Tits metric $|\cdot,\cdot|_T$. The Tits boundary of a CAT(0) space is a CAT(1) space.
The Euclidean $n$-space has the round $(n-1)$-sphere as Tits boundary, while any hyperbolic space has a discrete Tits boundary.
 The {\em Tits cone} of $X$ is denoted by $C_T(X)$,
it is the Euclidean cone over $\tits X$. For any point $x\in X$ there are  natural 1-Lipschitz exponential and logarithm maps
\[\exp_x:C_T(X)\to X\quad\quad \log_x:\tits X\to\Si_x X.\]
They satisfy the following rigidity. If $\log_x$ is isometric on a subset $A\subset\tits X$, then $\exp_x$
restricts to an isometric embedding on $C_T(A)\subset C_T(X)$ \cite[Flat Sector Lemma~2.3.4]{KleinerLeeb}.
 
A subset in a CAT(0) space is convex, if it contains the geodesic between any pair of its points.

If $C$ is a closed convex subset, then it is CAT(0) with respect  to the induced metric. In this case, $C$ admits a 1-Lipschitz retraction $\pi_C:X\to C$.
If $C_1$ and $C_2$ are closed convex subsets, then the distance function $d(\cdot,C_1)|_{C_2}$ is convex; and constant if
and only if $\pi_{C_1}$ restricts to an isometric embedding on $C_2$. We call $C_1$ and $C_2$ {\em parallel}, $C_1\| C_2$, if
and only if  $d(\cdot,C_1)|_{C_2}$ and $d(\cdot,C_2)|_{C_1}$ are constant.
Let $Y\subset X$ be a geodesically complete closed convex subset. Then we define the {\em parallel set} $P(Y)$
as the union of all closed convex subsets parallel to $Y$. The parallel set is closed, convex and splits canonically as a direct product
\[P(Y)\cong Y\times CS(Y)\]
where the {\em cross section} $CS(Y)$ is a closed convex subset.
The Tits boundary is given by
\[\tits P(Y)\cong\tits Y\circ\pol(\tits Y).\]
If $X_1$ and $X_2$ are CAT(0) spaces, then their direct product $X_1\times X_2$ is again a CAT(0) space.
We have $\tits (X_1\times X_2)=\tits X_1\circ \tits X_2$ and $\Si_{(x_1,x_2)}(X_1\times X_2)=\Si_{x_1} X_1\circ \Si_{x_2} X_2$.  
If $X$ is a geodesically complete CAT(0) space, then any join decomposition of $\tits X$ is induced by a direct product decomposition of $X$ \cite[Proposition~2.3.7]{KleinerLeeb}.
A CAT(0) space $X$ is called {\em irreducible}, if it does not admit a non-trivial splitting  as a direct product.

Let $X_1$ and $X_2$ be CAT(0) spaces containing closed convex subsets $C_1$ and $C_2$, respectively.
If there is an isometry $f:C_1\to C_2$, then by \cite[8.9.1 Reshetnyak's gluing theorem]{AKP}, the glued space
\[X_1\cup_f X_2:=(X_1\dot\cup X_2)/x_1\sim f(x_1)\]
is CAT(0) with respect to the induced length metric.
In particular, for every CAT(0) space $X$ which contains a closed convex subset $C$ we can construct a CAT(0) space $\hat X$, the {\em double}
of $X$ along $C$, by
\[\hat X=X^-\cup_C X^+\]
where $X^\pm$ are two isometric copies of $X$. 

The Tits boundary of a closed convex subset $Y\subset X$ embeds canonically $\tits Y\subset\tits X$.
If two closed convex subsets $Y_1$ and $Y_2$ intersect in $X$, then 
\[\tits(Y_1\cap Y_2)=\tits Y_1\cap \tits Y_2.\]

For points $\xi\in\geo X$ and $x\in X$ we denote the {\em Busemann function centered at $\xi$ based at $x$} by $b_{\xi,x}$.
If $\rho:[0,\infty)\to X$ denotes the geodesic ray asymptotic to $\xi$ and with $\rho(0)=x$, then
\[b_{\xi,x}(y)=\lim\limits_{t\to\infty}(|y,\rho(t)|-t).\]
It is a 1-Lipschitz convex function whose negative gradient at a point $y\in X$ is given by $\log_y(\xi)$. 
We denote the {\em horoball} centered at a point $\xi\in\geo X$ and based at the point $x\in X$ by 
\[HB_{\xi,x}:=b^{-1}_{\xi,x}([0,\infty)).\]
It is a closed convex subset with 
\[\tits HB_\xi(x)=\bar B_{\frac{\pi}{2}}(\xi)\subset\tits X.\]

A {\em $n$-flat} $F$ in a CAT(0) space $X$ is a closed convex subset isometric to $\R^n$.
In particular, $\tits F\subset\tits X$ is a round $(n-1)$-sphere.
On the other hand, if $X$ is locally compact and $\si\subset \tits X$ is a round  $(n-1)$-sphere, then either there exists an
$n$-flat $F\subset X$ with $\tits F=\si$, or there exists a round $n$-hemisphere $\tau^+\subset\tits X$
with $\si=\D\tau^+$ \cite[Proposition~2.1]{Leeb}. Consequently, if $\tits X$ is $(n-1)$-dimensional, then any round $(n-1)$-sphere in $\tits X$
is the Tits boundary of some $n$-flat in $X$. Moreover, a round $n$-hemisphere $\tau^+\subset\tits X$ bounds a flat $(n+1)$-half-space in $X$ if and only if
its boundary $\D\tau^+$ bounds an $n$-flat in $X$. A {\em flat ($n$-dimensional) half-space} $H\subset X$ is a closed convex subset isometric to a Euclidean half-space $\R_+^n$.
Its boundary $\D H\subset X$ is an $(n-1)$-flat and its Tits boundary is a round $(n-1)$-hemisphere $\tits H\subset\tits X$.  Flat half-spaces will play a certain role in our arguments later, and we agree to denote them by $H$
and there boundaries by $\D H=h$.

We define the {\em parallel set} of a round sphere $\si\subset \tits X$ as 
\[P(\si)=P(F)\]
where $F$ is a flat in $X$ with $\tits F=\si$, if such a flat exists. 

The {\em rank} of a CAT(0) space $X$ is defined by\footnote{There is no standard definition of rank for singular spaces in the literature,
several suggestions were made in \cite{G_asym} and their relation has been investigated in \cite{Kleiner}.} 
\[\rank(X)=\dim(\tits X)+1.\]
If $X$ is a locally compact CAT(0) space with an isometry group that acts cocompactly,
then the rank of $X$ coincides with the maximal dimension of a flat in $X$ \cite[Theorem~C]{Kleiner}.

\section{Buildings}\label{sec_buildings}

\subsection{Metric characterizations}
Originally, buildings were introduced by Tits in order to provide geometric interpretations for algebraic groups.
A geometric treatment of spherical and Euclidean buildings from the perspective of CAT($\ka$)-spaces
was developed in \cite{KleinerLeeb}. Here we take the latter point of view.
The relevance of spherical and Euclidean buildings to Alexandrov geometry lies in the fact that they appear as  Tits boundaries
and asymptotic cones, respectively, of non-positively curved symmetric spaces of higher rank \cite{KleinerLeeb}.
Instead of recalling the usual definitions of spherical and Euclidean buildings (which can be found in \cite{KleinerLeeb})
we provide simple metric characterizations proved in \cite{CL_metric}.

\bdfn
Let $X$ be a connected, piecewise spherical (respectively Euclidean)
complex of dimension $n\geq 2$ satisfying
\begin{enumerate}
	\item  $X$ is CAT(1) (respectively CAT(0)).
\item Every $(n-1)$–cell is contained in at least two $n$-cells.
\item Links of positive dimension are connected.
\item Links of dimension 1 have diameter $\pi$.
\end{enumerate}
Then $X$ is isometric to a spherical building (respectively a Euclidean
building).
\edfn

\bdfn
A {\em 1-dimensional spherical building} is a 1-dimensional piecewise spherical CAT(1) space of diameter equal to $\pi$
where every vertex has valence at least 2.
A {\em 1-dimensional Euclidean building} is a geodesically complete tree.
\edfn

Note that every building is geodesically complete. Moreover, in an $n$-dimensional spherical building (respectively Euclidean
building) any pair of points is contained in a round $n$-sphere (respectively an $n$-flat).
In particular, spherical buildings have diameter $\pi$ and the rank of a Euclidean building coincides with its dimension. 
However, a piecewise spherical CAT(1) spaces of diameter $\pi$ does not have to be a spherical building \cite{BB_diam}.
Similarly, a geodesically complete, piecewise Euclidean CAT(0) complex where every pair of points lies in a 2-flat
does not have to be a product or a Euclidean building \cite{BB_rr}.  

Examples of spherical buildings arise from projective planes, whereas
examples of (locally compact irreducible) Euclidean buildings 
arise from simple algebraic groups 
over non-Archimedean locally compact fields with a discrete valuation, see for instance \cite[Section~2.3]{Leeb}.

\subsection{Recognizing a building at infinity}

The aim of this section is to prove the following result which will be used in the proof of Theorem~\ref{thm_mainA}.

\bthm\label{thm_rec}
Let $X$ be a locally compact, geodesically complete CAT(0) space of rank $n$. Suppose that every complete geodesic in $X$ lies in an $n$-flat.
It $\tits X$ contains an open
relatively compact subset $U$, then $\tits X$ is a building.
\ethm

The result is similar to \cite[Theorem~1.6]{BL_building} and builds on Theorem~\ref{thm_lytchak} in the same way.
Note that if we would know that any pair of antipodes in $\tits X$ lie in a round $(n-1)$-sphere, 
then we could directly apply \cite[Theorem~1.6]{BL_building}

We need some preparation before we can provide the proof of Theorem~\ref{thm_rec}.

\blem\label{lem_ext_angle}
Let $X$ be a locally compact, geodesically complete CAT(0) space. Let $c$ be a complete geodesic in $X$ and let $\xi$ be a point in $\geo X\setminus\D c$.
Then for every $p\in c$ there exists a complete geodesic $c'$ extending the geodesic ray $p\xi$ and such that 
\[\angle_p(c(-\infty),\hat\xi)=\pi-\angle_p(c(-\infty),\xi)\]
where $\hat\xi:=c'(+\infty)$.
\elem

\proof
Since $X$ is locally compact and geodesically complete, all links in $X$ are geodesically complete \cite[Corollary~5.9]{LN_gcba}.
Denote by $v^\pm\in \Si_p X$ the direction at $p$ pointing to $c(\pm\infty)$, and let $w\in \Si_p X$ be the direction at $p$
pointing to $\xi$. We can extend the geodesic from $w$ to $v^-$ to an antipode $\hat w$ of $w$. Again, since $X$ is locally compact and geodesically complete,
there exists a geodesic ray starting in $p$ in the direction $w$. This proves the claim.
\qed
\medskip

Recall that a pair of antipodes in the Tits boundary of a CAT(0) space $X$ does not have to bound a complete geodesic in $X$.
Nevertheless, for locally compact geodesically complete spaces we have the following.

\blem\label{lem_antipodes}
Let $X$ be a locally compact, geodesically complete CAT(0) space. Let $\xi$ and $\hat\xi$ be antipodes in $\tits X$.
Then there exist a sequence of complete geodesics $c_k$ in $X$ with $\geo c_k=\{\hat\xi,\xi_k\}$ and $\xi_k\to\xi$ in $\geo X$.
\elem

\proof
Fix a base point $o\in X$ and let $\rho:[0,\infty)\to X$ be a parametrization of the geodesic ray $o\xi$.
By Lemma~\ref{lem_ext_angle}, we can choose for each $k\in\N$ a complete geodesic $c_k$ extending the geodesic ray $\rho(k)\hat\xi$ and such that 
\[\angle_{\rho(k)}(\hat\xi,\xi)=\pi-\angle_{\rho(k)}(\xi,\xi_k)\]
where $\xi_k:=c_k(+\infty)$. 
Now for every $p\in o\xi$ we have
\[\angle_p(\xi_k,\xi)\leq\angle_{\rho(k)}(\xi_k,\xi)=\pi-\angle_{\rho(k)}(\xi,\hat\xi)\to 0.\]
Since $\angle_{\rho(k)}(\xi,\hat\xi)\to|\xi,\hat\xi|=\pi$. As $p\in o\xi$ was arbitrary, this shows $\xi_k\to\xi$ in $\geo X$.
\qed
\medskip

\blem\label{lem_anti_spheres}
Let $X$ be a locally compact, geodesically complete CAT(0) space of rank $n$. Suppose that every complete geodesic in $X$ lies in an $n$-flat.
If $\xi\in\tits X$ is a regular point, then for every antipode $\hat\xi$ of $\xi$ there exists a round $(n-1)$-sphere $\si\subset\tits X$
containing $\xi$ and $\hat\xi$.
\elem

\proof
Let
$\xi\in\tits X$ be a regular point. and $\si\subset\tits X$ a $(n-1)$-sphere containing $\xi$.
Then there exists $s>0$ such that $B_s(\xi)$ is contained in any $(n-1)$-sphere which contains $\xi$.
If $\hat\xi\in\tits X$ is an antipode of $\xi$, then by Lemma~\ref{lem_antipodes}, there exists a sequence of complete geodesics $(c_k)$
in $X$ with $\geo c_k=\{\xi,\hat\xi_k\}$ and $\hat\xi_k\to\hat\xi$. By assumption, there exists $n$-flats $F_k$ with $c_k\subset F_k$. 
Hence,  we obtain $|\eta,\hat\xi|=\pi-s$
for every $k\in\N$ and $\eta\in\tits X$ with $|\eta,\xi|=s$. Thus the claim follows from Lemma~\ref{lem_lune}.
\qed

\bprop\label{prop_reg_anti}
Let $X$ be a locally compact, geodesically complete CAT(0) space of rank $n$. Suppose that every complete geodesic in $X$ lies in an $n$-flat.
Then the set $O\subset \tits X$ of regular points is closed under taking antipodes. 
\eprop

\proof
Let $\xi$ be a regular point in $\tits X$ and let $\hat \xi$ be an antipode.
By Lemma~\ref{lem_anti_spheres}, there exists a round $(n-1)$-sphere $\si\subset \tits X$ which contains $\xi$
and $\hat\xi$. Since $\xi$ is regular, there exists $s>0$ such that $B_s(\xi)\subset\si$.
We claim that $B_s(\hat\xi)\subset\si$ holds as well.
Let $\hat\eta$ be a point in $B_s(\hat\xi)$.
By Lemma~\ref{lem_gc}, $\tits X$ is geodesically complete and we can extend the geodesic $\hat\eta\hat\xi$
up to an antipode $\eta$ of $\xi$ in $\si$. Then $|\eta,\xi|=|\hat\eta,\hat\xi|<s$ and therefore $\eta$ is regular.
By Lemma~\ref{lem_anti_spheres}, we find a round $(n-1)$-sphere $\si'\subset \tits X$ which contains $\eta$
and $\hat\eta$. Because $\eta\in B_s(\xi)$, we must have $B_s(\xi)\subset\si\cap\si'$. By construction,
$\hat \xi$ also lies in $\si\cap\si'$ and thus $\si=\si'$ and $\hat\eta\in\si$ as required.
\qed

\blem\label{lem_split}
Let $X$ be a locally compact, geodesically complete CAT(0) space of rank $n$. 
Suppose $X$ splits as a non-trivial direct product $X\cong X_1\times X_2$.
Then for $j=1,2$, we have $n=n_1+n_2$ where $n_j\geq 1$ denotes the rank of $X_j$, 
and every complete geodesic in $X_j$ is contained in an $n_j$-flat.
\elem

\proof
It is clear that the rank is additive. So
let $c_j$ be a complete geodesic in $X_j$ and let $\tilde c=(c_1,c_2)$ be the diagonal geodesic in
$c_1\times c_2\subset X_1\times X_2$. If $\tilde F$ is an $n$-flat in $X$ which contains $\tilde c$,
then by \cite[Lemma~2.3.8]{KleinerLeeb} there are flats $F_j\subset X_j$ such that $\tilde F\subset F_1\times F_2$.
In particular, $\dim(F_j)=n_j$ and $c_j\subset F_j$. 
\qed
\medskip

\proof[Proof of Theorem~\ref{thm_rec}]
By Lemma~\ref{lem_gc}, $\tits X$ is geodesically complete and has diameter $\pi$.

We proceed by induction on the rank of $X$. If the rank is 1, then $\tits X$ is discrete and there is nothing to show. 
Suppose $X$ splits as a non-trivial direct product $X\cong X_1\times X_2$.
Then by Lemma~\ref{lem_split},  the factors  are subject to the induction
hypothesis and we may assume that
$X$ is irreducible.

If $U$ is an open relatively compact subset of $\tits X$, then $U$ contains a subset $U'$
homeomorphic to an $n$-manifold \cite[Theorem~1.2]{LN_gcba}. Hence, for a point $\xi\in U'$
there exists an $\eps>0$ such that $B_\eps(\xi)\subset U'$. If $\si\subset \tits X$
denote a round $(n-1)$-sphere with $\xi\in\si$, then $\si\cap B_\eps(\xi)=B_\eps(\xi)$.
In particular, $\xi$ is regular and the regular set  $O\subset\tits X$
is non-empty. If $O=\tits X$, then $X$ is isometric to $\R^n$. So let us assume that $O$ is a proper
subset of $\tits X$. By Proposition~\ref{prop_reg_anti}, the set $\tits X\setminus O$ is
a non-empty, closed proper subset which contains all antipodes of all of its points.
Since $\tits X$ is irreducible, Theorem~\ref{thm_lytchak} implies that it is a building.
\qed

\section{Orthogonal flat half-planes vs. local compactness}\label{sec_illu}

This short section is illustrative. We present the key observation behind discreteness of branch points in a simplified setting.

The following basic fact will be used repeatedly throughout the paper.

\blem[{\cite[Sublemma~2.3]{Leeb}}]\label{lem_flats_in_convex}
Let $F$ be a flat in a CAT(0) space and let $C\subset X$ be a closed convex subset so that $\geo F\subset\geo C$. 
Then $C$ contains a flat $F'$ parallel to $F$.
\elem

\proof
The function $d_C$, the distance function to $C$, is a convex function on $F$. By assumption, $d_C$ is bounded on $F$ and therefore constant, $d_C\equiv a$.
By \cite[Flat Strip Lemma~2.3.5]{KleinerLeeb}, there exists an isometric embedding $\psi:F\times[0,a]\to X$
such that $\psi(\cdot,0)=\id_F$ and $\psi(\cdot,a)=\pi_C|_{F}$. This implies the claim.  
\qed

\medskip
In the upcoming section we will see how branch points at infinity of a flat $F$ lead to flat half-planes orthogonal to $F$. 
Hence the next lemma, which is not needed in the rest of the paper, explains the conflict between branch points at infinity and local compactness of a space. 
The responsible geometric property is flatness of intersecting parallel sets, similar as in \cite[Theorem~2]{St_obst}.

\blem\label{lem_illu}
Let $X$ be a CAT(0) space of rank 2 and let $F\subset X$ be a 2-flat with $P(F)=F$.
Suppose that $H^\pm\subset X$ are two flat half-planes, orthogonal to $F$.
If the boundaries $\D H^\pm\subset F$ are not parallel, then the angle between $H^-$
and $H^+$ is at least $\frac{\pi}{2}$. 
\elem

\begin{center}
\includegraphics[scale=0.5,trim={0cm 0cm 0cm 0cm},clip]{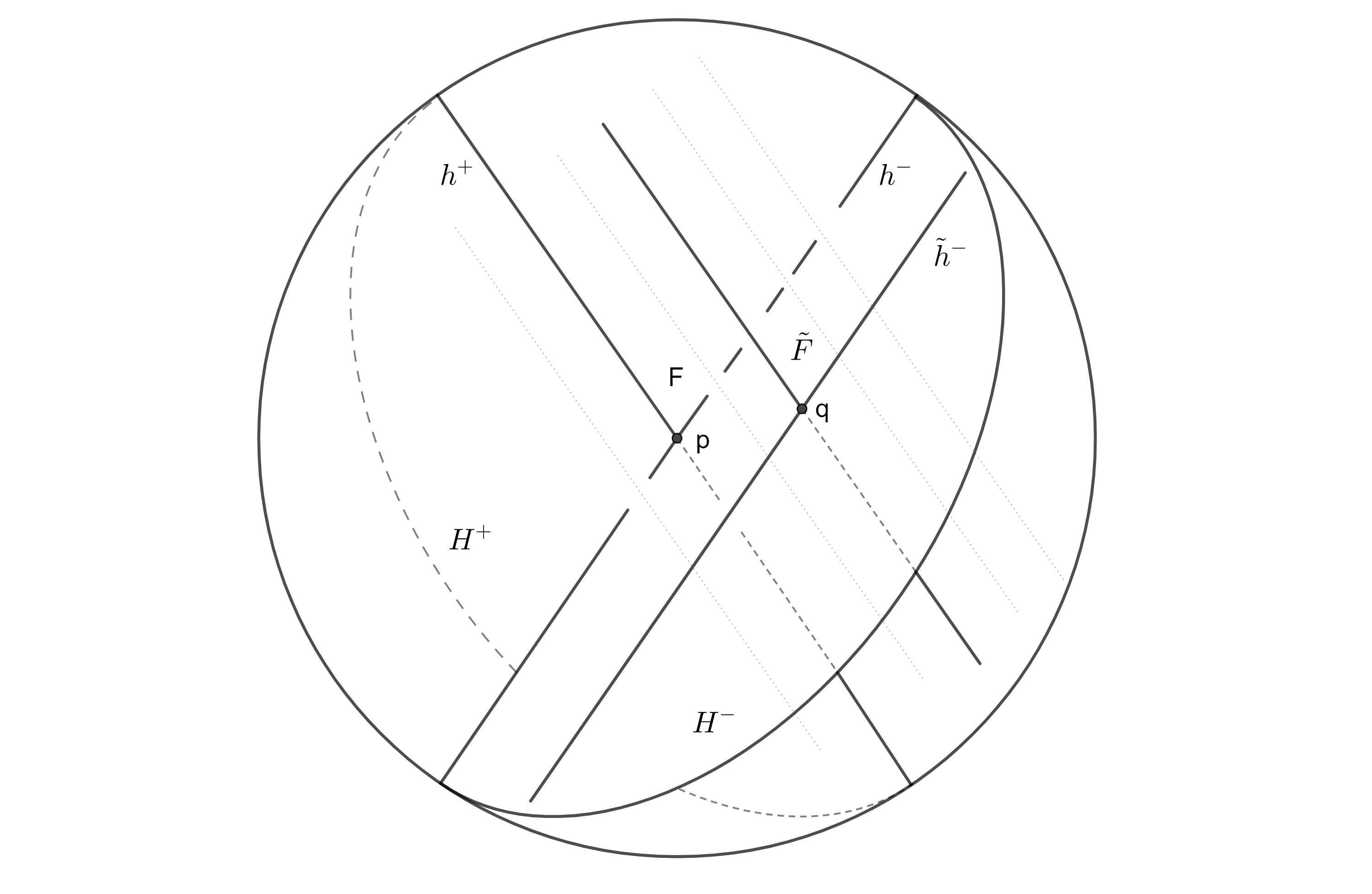}
\end{center}

\proof
Set $h^\pm:=\D H^\pm$ and put $p=h^-\cap h^+$. Further, denote by $\eta^\pm$ the centers of $\geo H^\pm$.
Since $H^-$ is orthogonal to $F$, the Busemann function $b_{\eta^+,p}$ is non-negative on $F$.
If $\angle_p(\eta^-,\eta^+)<\frac{\pi}{2}$, then, by the first variation formula, there exists a point $q\in H^-$
such that $b_{\eta^+,p}(q)<0$. Since $H^-\subset P(h^-)$, we see $HB(\eta^+,q)\cap P(h^-)\neq\emptyset$.
Because $\geo h^+\subset \geo (HB(\eta^+,q))\cap \geo P(h^-))=\geo (HB(\eta^+,q)\cap P(h^-))$,
we infer from Lemma~\ref{lem_flats_in_convex} that there is a complete geodesic $\tilde h^-\subset HB(\eta^+,q)\cap P(h^-)$ which is parallel to $h^-$.
Since $h^-$ is not parallel to $h^+$, there exists a 2-flat $\tilde F\subset P(h^-)$ which contains $\tilde h^-$.
Thus, $\geo\tilde F=\geo F$ and therefore $\tilde F$ is parallel to $F$. By assumption, we must have $\tilde F=F$.
This is a contradiction because  $b_{\eta^+,p}$ is bounded above on $HB(\eta^+,q)$ by its value at $q$. In particular, it is negative on $\tilde h^-$.
But at the same time it is non-negative on $F$.
\qed

\section{Rank $n$ spaces with periodic $n$-flats}\label{sec_rank_n}

\subsection{From branch points to orthogonal half-planes}

\bdfn
Let $X$ be a locally compact CAT(0) space and $F\subset X$ a flat.
Let $\Isom(X)$ denote the isometry group of $X$.
We say that $F$ is {\em periodic}, if its {\em stabilizer} $\Stab(F):=\{\ga\in \Isom(X)|\ \ga(F)=F\}$ contains a subgroup acting geometrically on $F$.
\edfn

\blem\label{lem_bieber}
Let $X$ be a locally compact CAT(0) space and let $F\subset X$ be a periodic $n$-flat with $n\geq 2$.
Suppose that $F=F^-\times F^+$ is a direct product decomposition into flats $F^\pm$ of dimension $n^\pm$.
Then $F^\pm$ is a periodic $n^\pm$-flat in the cross section $CS(F^\mp)$.
\elem

\proof
By Bieberbach's theorem, the stabilizer of $F$ contains a free abelian subgroup $A$ of rank $n$
which acts cocompactly on $F$ by Euclidean translations. Since $A$ preserves $\tits F$ pointwise, it also preserves $P(F^\mp)$
and the splitting $P(F^\mp)\cong F^\mp\times CS(F^\mp)$.
Hence its projection to $CS(F^\mp)$ contains a free abelian subgroup of rank $n^\pm=n-n^\mp$ which  acts 
geometrically by Euclidean translations on $F^\pm\subset CS(F^\mp)$. 
\qed
\medskip

The following result implies that a branch point in the ideal boundary of a periodic flat leads to an orthogonal flat half-plane.

\blem[{\cite[Lemma~2.3.1]{HK}}]\label{lem_perpH}
Let $X$ be a locally compact CAT(0) space.
Let $F\subset X$ be a periodic flat. Let $\xi\in\tits F$ be a point such that there exists a sequence $(\xi_k)$ in $\tits X\setminus\tits F$
which converges to $\xi$ with respect to the Tits metric. 
Then there exists a flat half-plane $H\subset X$ asymptotic to $\xi$ and orthogonal to $F$. 
More precisely, $H$ satisfies $\D H\subset F$, $\xi\in\geo(\D H)$ and
$\angle(H,F)\geq \frac{\pi}{2}$. 
\elem

\proof
Since $F$ is periodic, it is enough to find for every $R>0$ a flat strip $S\cong\R\times[0,R]$ which is asymptotic to $\xi$ and
orthogonal to $F$. We fix a base point $o\in F$. For every $k\in\N$ we choose $x_k\in o\xi_k$ such that $|x_k,F|=R$.
Denote by $\bar x_k\in F$ the nearest point to $x_k$. Choose isometries $\ga_k\in\Stab(F)$ such that all points $\ga_k(\bar x_k)$
lies in a fixed compact set. After passing to a subsequence, we obtain a limit geodesic $c_\infty$ which is parallel to $F$
and spans the required flats strip. More precisely, we have the convergence $\ga_k(o\xi_k,x_k)\to(c_\infty,x_\infty)$ with respect to pointed Hausdorff topology.   
The complete geodesic $c_\infty$ is asymptotic to $\xi$ since $|\ga_k\xi_k,\xi|=|\xi_k,\xi|\to 0$. Because $x_k o\subset \bar N_R(F)$, also $c_\infty$
is contained in the $R$-neighborhood of $F$ and has even distance $R$ from $F$ since $|x_\infty,F|=R$.
Denote by $\bar c_\infty\subset F$ the complete geodesic which contains the point $x_\infty$ and is asymptotic to $\xi$.
Then $c_\infty$ and $\bar c_\infty$ span the required flat strip by \cite[Lemma~2.3.5]{KleinerLeeb}.  
\qed

\subsection{Orthogonal half-spaces from a sequence of orthogonal half-planes}

\bprop\label{prop_halfspace}
Let $X$ be a locally compact  CAT(0) space.
Let $\hat F$ be a periodic $\hat n$-flat in $X$ with $\tits\hat F=\hat \si$.
Further, let $\si\subset\hat\si$ be a round 1-sphere. Suppose that there is a sequence of pairwise non-parallel flat half-planes $(H_k)$,
orthogonal to $\hat F$, with boundary $h_k:=\D H_k$ and with $\geo h_k \subset\si$. Then there exists a 3-dimensional flat half-space $\hat H$, orthogonal to $\hat F$,
with boundary $\hat h:=\D \hat H$ and with $\geo \hat h=\si$.
\eprop

\proof
Since $\hat F$ is periodic, we may assume that we have a convergent sequence of pairwise distinct pointed flat half-planes $(H_k,p_k)\to (H_\infty,p_\infty)$, 
where $p_k\in h_k\cap\hat F$.
The limit $H_\infty$ with boundary $h_\infty:=\D H_\infty$ is then still orthogonal to $\hat F$ and $\geo h_\infty\subset\si$. 
Set $\tits H_k=\tau_k^+$, $\tits H_\infty=\tau_\infty^+$ and denote by $\zeta_k$ and $\zeta_\infty$
the respective centers. Because $p_\infty\in h_\infty$ and $H_\infty$ is orthogonal to $\hat F$, the Busemann function $b_{\zeta_\infty,p_\infty}$
is non-negative on $\hat F$. Since the geodesic rays $p_k\zeta_k$ converge to $p_\infty\zeta_\infty$, we find a sequence $t_k\to\infty$ and points $x_k\in p_k\zeta_k$
such that $b_{\zeta_\infty,p_\infty}(x_k)\leq-t_k$. The geodesic ray $p_k\zeta_k$ lies in $H_k$ and therefore $x_k$ lies in $P(h_k)$.
Thus $HB(\zeta_\infty,x_k)\cap P(h_k)$ is a non-empty closed convex set and we have
\[\geo(HB(\zeta_\infty,x_k)\cap P(h_k))=\bar B_{\frac{\pi}{2}}(\zeta_\infty)\cap\geo P(h_k).\]
As $\hat F\subset P(h_k)$, we see $\geo h_\infty\subset \geo(HB(\zeta_\infty,x_k)\cap P(h_k))$.
From Lemma~\ref{lem_flats_in_convex} we infer that there is a complete geodesic $c_k\subset HB(\zeta_\infty,x_k)\cap P(h_k)$
with $\geo c_k=\geo h_\infty$. Since the flat half-planes $H_k$ are pairwise distinct and $c_k\subset P(h_k)$,
we see that there is a 2-flat $F_k$ in $P(h_k)$ which contains $c_k$ and such that $\geo F_k=\si$.
Note  that $b_{\zeta_\infty,p_\infty}$ is bounded above on $c_k$ by $-t_k$ since $c_k\subset HB(\zeta_\infty,x_k)$.
As $b_{\zeta_\infty,p_\infty}$ is non-negative on $\hat F$ we conclude $c_k\cap N_{t_k}(\hat F)=\emptyset$. Hence if $F'_k\subset\hat F$
denotes the closest parallel flat to $F_k$, then $|F_k,F'_k|_H\geq t_k$. In particular, $F_k$ and $F'_k$
bound a flat strip orthogonal to $\hat F$ and of width at least $t_k$. Since $\hat F$ is periodic and $t_k\to\infty$, we obtain the claim.   
\qed
\medskip

An immediate consequence is the following  finiteness result for spaces of rank 2.

\bcor\label{cor_finite_branch}
Let $X$ be a locally compact CAT(0) space of rank 2.
Let $F\subset X$ be a periodic 2-flat. Then $\si:=\tits F\subset\tits X$ contains only a finite set  of branch points $E$
its complement $\si\setminus E$ is open in $\tits X$.
\ecor

\proof
By Lemma~\ref{lem_perpH}, every branch point $\xi$ in $\si$ yields a flat half-plane orthogonal to $F$ and
asymptotic to $\xi$. Since $X$ has rank 2, Proposition~\ref{prop_halfspace} shows that there are only finitely many such half-planes. 
\qed

\subsection{Periodic flats lead to regular points at infinity}

\blem\label{lem_regline}
Let $X$ be a locally compact CAT(0) space and let $F\subset X$ be an $n$-flat, $n\geq 2$, with $\tits F=\si$.
If $\si$ contains a non-empty open subset $U\subset\tits X$, then any complete geodesic $l\subset F$
asymptotic to $U$ satisfies $\tits P(l)=\si$. In particular, $P(l)$ lies in a finite tubular neighborhood of $F$.
\elem

\proof
We write $F=l\times F'$.
Recall that $\tits P(l)\cong\tits l\circ\tits CS(l)$. Hence if there is a point $\eta\in\tits CS(l)\setminus\tits F'$,
then no neighborhood of $l(+\infty)$ in $\si$ can be open in $\tits X$.
Contradiction.  
\qed

\bthm\label{thm_open}
Let $X$ be a locally compact CAT(0) space of rank $n\geq 2$.  If $X$ contains a periodic $n$-flat $\hat F$, then $\hat\si=\tits \hat F$ contains a dense subset
$U\subset \hat\si$ which is open in $\tits X$.
\ethm

\proof
We will prove the statement by induction on the rank of $X$.
The case of rank 2 follows from Corollary~\ref{cor_finite_branch}.
Now suppose the result holds for all locally compact CAT(0) spaces of rank at most $n-1$.

Choose a point $\xi^+$ in $\hat\si$.
Denote by $\xi^-$ its antipode in $\hat\si$ and write $\hat\si=\{\xi^-,\xi^+\}\circ\si'$.
Let $\hat F=l\times F'$ be a corresponding splitting. Note that the cross section $CS(\{\xi^-,\xi^+\})$
is a locally compact CAT(0) space of rank $(n-1)$ which contains the periodic $(n-1)$-flat $F'$ (Lemma~\ref{lem_bieber}).
By induction hypothesis and Lemma~\ref{lem_regline}, there exists a complete geodesic $l'\subset F'$ such that its parallel set in
$CS(\{\xi^-,\xi^+\})$ lies in a finite tubular neighborhood of $F'$. Define $F:=l\times l'\subset \hat F$ and set $\si:=\tits F\subset\hat\si$.
Then the parallel set $P(\si)$ lies in a finite tubular neighborhood of $\hat F$. In particular, there is no 3-dimensional flat half-space $\hat H$,
orthogonal to $\hat F$, with boundary $\hat h:=\D\hat H$ and $\geo \hat h=\si$. Hence by Proposition~\ref{prop_halfspace},
there exists a finite set $E\subset\si$ such that if a pair of antipodes $\{\xi_0^-,\xi_0^+\}\subset \si$
is not contained in $E$, then the parallel set $P(\{\xi_0^-,\xi_0^+\})$ lies in a finite tubular neighborhood of $\hat F$.

Now we choose a pair of antipodes $\{\xi_0^-,\xi_0^+\}\subset \si$ which is not contained in $E$.
Then a small open neighborhood of $\xi_0^+$ in $\hat\si$ has to be open in $\tits X$.
Otherwise, there exists a sequence of branch points $\xi_k^+$ in $\hat\si$ with $\xi_k^+\to\xi_0^+$.
By Lemma~\ref{lem_perpH}, and since $\hat F$ is periodic, there exists a convergent sequence $(H_k)$ of flat half-planes, orthogonal to $\hat F$, 
and with boundaries $h_k:=\D H_k$ and $\xi_k^+\in\geo h_k$. Their limit $H_\infty$ is a flat half-plane, orthogonal to $\hat F$, with boundary $h_\infty:=\D H_\infty$  
and $\xi_0^+\in\geo h_\infty$. In particular, $H_\infty\subset P(\{\xi_0^-,\xi_0^+\})$. Contradiction.

It follows that all points in $\si\setminus E$ are regular.
Thus, $\xi^+$ lies in the closure of the regular points in $\si$. In particular,
since $\xi^+\in\hat \si$ was arbitrary, this shows that $\hat \si$ intersects the set of regular points in a dense subset.
\qed
\medskip

\proof[Proof of Theorem~\ref{thm_mainA}]
Let $F\subset X$ be a periodic $n$-flat.
By Theorem~\ref{thm_open}, $\si=\tits F$ contains a  dense subset which is open in $\tits X$.
Then, by Theorem~\ref{thm_rec}, $\tits X$ is a spherical join or a spherical building.
Since $X$ is geodesically complete and locally compact, the claim follows from Theorem~\ref{thm_leeb} 
and \cite[Proposition~2.3.7]{KleinerLeeb}.
\qed

\proof[Proof of Corollary~\ref{cor_mainA}]
By Theorem~\ref{thm_mainA}, Lemma~\ref{lem_split} and Lemma~\ref{lem_bieber}, we can produce a product decomposition of $X$ as claimed.
Since every such decomposition induces a corresponding join decomposition of $\tits X$.
The uniqueness statement follows from the uniqueness of such join decompositions \cite[Corollary~1.2]{Ly_rigidity}.
\qed

\section{Rigidity in rank 2}\label{sec_rank_2}

\subsection{Projecting ideal points onto parallel sets}

\blem\label{lem_minimum}
Let $X$ be a locally compact CAT(0) space of rank 2. 
Let $\xi^-$ and $\xi^+$
be a pair of antipodes in $\tits X$ and let $\eta$ be a pole of $\{\xi^-,\xi^+\}$.
Further, let $P\subset X$ be a closed convex subset with $\eta\notin\tits P$ and $\{\xi^-,\xi^+\}\subset\tits P$.
If $b_\eta$ is a horofunction associated to $\eta$, then $b_\eta$ attains a minimum on $P$. Moreover, the minimum set 
contains a complete geodesic $c$ with $\geo c=\{\xi^+,\xi^-\}$.  
\elem

\proof
Let us first show that the parallel set $P(\xi^+,\xi^-)$ is non-empty. Denote by $\tau^+\subset \tits X$
the round hemisphere with $\D\tau^+=\{\xi^+,\xi^-\}$ and center $\eta$.
If there is no round hemisphere $\tilde\tau^+$ in $\tits P$ with $\D\tilde\tau^+=\{\xi^+,\xi^-\}$, then by \cite[Proposition~2.1]{Leeb},
there exists a complete geodesic $c$ in $P$ with $\geo c=\{\xi^+,\xi^-\}$. On the other hand, if there does exist such a 
$\tilde\tau^+$, then, since  $\eta\notin\tits P$, the union $\si:=\tilde\tau^+\cup\tau^+$ forms a round 1-sphere in $\tits X$.
Since $X$ has rank 2, \cite[Proposition~2.1]{Leeb} ensures the existence of a 2-flat $F$ with $\tits F=\si$ and therefore $F\subset P(\xi^+,\xi^-)$.
Now it follows from Lemma~\ref{lem_flats_in_convex}, if $HB_\eta$ is a horoball based at $\eta$ such that $HB_\eta\cap P$ is non-empty, then 
$HB_\eta\cap P$ contains a complete geodesic $c$  with $\geo c=\{\xi^+,\xi^-\}$.
Since $b_\eta$ is constant on such geodesics, in order to find the required minimum, it is enough to restrict $b_\eta$
to $P\cap CS(\xi^+,\xi^-)$. Note that $|\eta',\eta|=\pi$ for every point $\eta'\in\tits CS(\xi^+,\xi^-)\setminus\{\eta\}$. 
Hence $\lim\limits_{x\to\eta'}b_\eta(x)=+\infty$ and 
$b_\eta$ does attain a minimum at a point $p$ in $P$. 
By assumption, we have
\[\{\xi^-,\xi^+\}\subset\tits(HB(p,\eta)\cap P).\]
Thus, by Lemma~\ref{lem_flats_in_convex}, there is a  complete geodesic $c$ in $HB(p,\eta)\cap P$ with $\geo c=\{\xi^+,\xi^-\}$.
Since $b_\eta$ is bounded above on $HB(\eta,p)$ by its value at $p$  the claim follows.   
\qed

\blem\label{lem_proj}
Let $X$ be a locally compact CAT(0) space of rank 2. Let $\si\subset\tits X$ be a round 1-sphere.
Let $\tau^+\subset\tits X$ be a round hemisphere with $\tau^+\cap\si=\D\tau^+$.
Then there exists a 2-flat $F\subset X$ with $\tits F=\si$ and a flat half-plane $H$,  orthogonal to $F$ and with
$\tits H=\tau^+$. Moreover, for any geodesic $c\subset F$ which is not parallel to $\D H$, holds $\angle(H,P(c))\geq\frac{\pi}{2}$.
\elem

\proof
Since $X$ has rank 2, the parallel set $P(\si)$  is a non-empty closed convex subset of $X$
which splits isometrically as $P(\si)\cong\R^2\times C$ where $C$ is a compact CAT(0) space.
Let $\eta$ be the center of $\tau^+$ and let $b_\eta$ be an associated horofunction.
By Lemma~\ref{lem_minimum}, 
$b_\eta$ attains a minimum on $P(\si)$ and the minimum set $Z\subset P(\si)$ contains a complete geodesic $l$ with $\geo l=\D\tau^+$.
Let $H$ be the flat half-plane with $\tits H=\tau^+$ and $\D H=l$.
Because $l\subset P(\si)$, there exists a 2-flat $F\subset X$ with $\tits F=\si$ which contains $l$.
Since $l$ lies in the minimum set $Z$, the half-plane $H$ has to be orthogonal to $F$.

Let $c\subset F$ be a geodesic not parallel to $l$. Choose  a point $x\in l$ and a geodesic ray $\rho$ starting in $x$ and asymptotic to $\eta$.
In particular, $\rho\subset H$.
 
We first claim that $\rho(t)\notin P(c)$ for $t>0$. Indeed, if $\rho(t)\in P(c)$, then, since $\rho(t)\in P(l)$, there exists lines $c_t\| c$ and $l_t\| l$ through
$\rho(t)$. These span a flat $F_t\subset P(\si)$. Hence $t=0$ because $\rho(0)$ minimizes $b_\eta$.

Now suppose $\angle(\rho,P(c))<\frac{\pi}{2}$. Then there exists a geodesic $c_1$ in $P(c)$ starting in $x$ and realizing the angle, 
$\angle(\rho,c_1)=\angle(\rho,P(c))$. In particular, $b_\eta(c_1(t))<b_\eta(x)$. 
Note that since $c$ is not parallel to $l$, the point $\eta$ is not in $\tits P(c)$. Moreover, we have $\D\tau^+\subset\tits P(c)\cap\bar B_{\frac{\pi}{2}}(\eta)$.
Hence, by Lemma~\ref{lem_minimum}, $b_\eta$ attains a minimum on $P(c)$, and the minimum set $Z_c$ contains a complete geodesic 
$\tilde l$ with $\D \tilde l=\{\xi, \hat\xi\}$.
But then $\tilde l$ lies in $P(c)$ and therefore, arguing as above, there is a 2-flat $\tilde F\subset P(\si)$ with $\tilde l\subset \tilde F$.
This contradicts the fact that $Z$ minimizes $b_\eta$ on $P(\si)$. Therefore $\angle(H,P(c))=\angle(\rho,P(c))\geq\frac{\pi}{2}$.
\qed

\subsection{Parallel sets do not accumulate in the Tits metric}\label{subsec_no_acc}

\blem\label{lem_no_acc}
Let $X$ be a locally compact CAT(0) space of rank 2. 
Let $g$ be  an axis of an axial isometry $\ga$. 
Suppose that there exists $a>0$ and a sequence of 2-flats $(F_k)$ in $P(g)$
with $g\subset N_a(F_k)$ for every $k\in\N$. Further, suppose that there is a sequence $(H_k)$
of flat half-planes with boundaries $h_k:=\D H_k\subset F_k$ and $\geo H_k\cap\geo F_k=\geo h_k$.
If $\geo h_k\to\geo g$ with respect to the Tits metric, then $\geo h_k=\geo g$
for almost all $k$.
\elem

\proof
Suppose for contradiction that the convergence $\geo h_k\to\geo g$ is non-trivial.
By Lemma~\ref{lem_proj}, there exist 2-flats $\tilde F_k\subset P(g)$ parallel to $F_k$; 
and flat half-planes $\tilde H_k$ orthogonal to $\tilde F_k$ and with $\geo \tilde H_k=\geo H_k$; and such that 
 $\tilde H_k$ is orthogonal to $P(g)$. Set $\tilde h_k:=\D \tilde H_k$. Since $X$ has rank 2 and a cocompact isometry group, there exists $D>0$
such that $|F_k,\tilde F_k|_H\leq D$ for all $k\in\N$. In particular, 
there exist complete geodesics $\tilde g_k\subset  \tilde F_k$ which are parallel to $g$ and
satisfy $|g, \tilde g_k|\leq a+D$ for all $k\in\N$. Denote by $\tilde x_k\in \tilde F_k$ the intersection point of $\tilde g_k$ and $ \tilde h_k$.
Then there are powers of $\ga$, denoted by $\ga_k$, such that the points $\ga_k(\tilde x_k)$ lie in a fixed compact set.
Since $\ga$ preserves $P(g)$, $\ga_k \tilde H_k$ is still orthogonal to $P(g)$. Since $X$ is locally compact, after choosing a subsequence, 
the flat half-planes $\ga_k \tilde H_k$ converge to a flat half-plane $\tilde H_\infty$. 
Note that by assumption $\ga_k (\geo \tilde h_k)\to\geo g$ with respect to the
Tits metric.
Hence $\tilde H_\infty\subset P(g)$. On the other hand, by upper semi-continuity of angles, $\tilde H_\infty$ has to be orthogonal to $P(g)$. Contradiction.
\qed

\subsection{Regular points at infinity}

In this section we will prove that the Tits boundary of a cocompact CAT(0) space of rank 2 contains regular points if its Tits boundary has diameter $\pi$.
To achieve this, we will use the following technical result in an essential way.

\begin{namedlemma}[Half-Plane]\label{lem_key_tech}
Let $X$ be a locally compact CAT(0) space of rank 2.  Let $\si\subset\tits X$ be a round sphere.
Further, let $g$ be an axis of an axial isometry $\ga$ and assume $\geo g\subset\si$.
Let $\xi^\pm\in\si$ be antipodes disjoint from $\geo g$.
Suppose that there is a sequence of local geodesics $\al^+_k\subset \tits X$ such that 
\begin{itemize}
	\item $\D\al^+_k=\{\xi^-_k,\xi^+_k\}$ and $\xi^\pm_k\to\xi^\pm$ with respect to the Tits metric;
	\item the lengths of $\al_k^+$  converges to $\pi$ as $k\to\infty$.
\end{itemize}
Then there exists a 2-flat $\tilde F\subset P(g)$, and a flat half-plane $\tilde H$ with boundary $\tilde h:=\D \tilde H\subset\tilde F$ such that the following properties hold.
\begin{enumerate}
	\item $g\subset N_a(\tilde F)$ if $g\subset N_a(F)$ for $a>0$;
	\item $|\geo \tilde h,\geo g|=|\{\xi^-,\xi^+\},\geo g|$;
	\item $\geo\tilde H\cap\geo\tilde F=\geo\tilde h$.
\end{enumerate}
\end{namedlemma}

The  proof of the \hyperref[lem_key_tech]{Half-Plane Lemma} requires methods from geometric measure theory.
Since these techniques do not play a role in the rest of the paper,
we defer their discussion, as well as the proof of the \hyperref[lem_key_tech]{Half-Plane Lemma}, to Appendix~\ref{sec_app}. 

\blem\label{lem_per2fl}
Let $X$ be a locally compact CAT(0) space with a geometric group action $\Ga\acts X$. 
Let $\ga\in\Ga$ be an axial isometry with axis $g$. Suppose that $g$ bounds a flat half-plane $H$
which is preserved by $\ga$. Then $X$ contains a periodic 2-flat $F\subset P(g)$.
\elem

\proof
Let $\rho\subset H$ be a geodesic ray orthogonal to $g$. Choose points $x_k\in\rho$ with $|x_k,g|\to\infty$.
We find a sequence $(\ga_k)\subset\Ga$ such that $\ga_k(x_k)$ stays in a fixed compact set $K$.
By the proof of \cite[Theorem~11]{Sw_cut}, we can pass to a subsequence such that for all $k<l$ the element $\ga_l^{-1}\ga_k$
is axial. Now note that since $\ga$ preserves $H$, its displacement  is constant on $H$, $|x,\ga x|=a>0$ for every $x\in H$.
Thus $|\ga_k x_k,(\ga_k\ga\ga_k^{-1})\ga_k x_k|=|x_k, \ga x_k|=a$. Since $\Ga$ is discrete and $\ga_k x_k\in K$, we may pass to a further subsequence
 such that $\ga_1\ga\ga_1^{-1}=\ga_k\ga\ga_k^{-1}$ holds for all $k\in \N$.
Hence $\beta_k=\ga_k^{-1}\ga_1$ is an axial isometry which commutes with $\ga$. The claim follows from \cite[Flat Torus Theorem~7.1]{BH}.
\qed

\bprop\label{prop_regular_points}
Let $X$ be a locally compact CAT(0) space of rank 2 with a geometric group action $\Ga\acts X$. 
Suppose that the diameter of $\tits X$ is equal to $\pi$.
Let $g$ be an axis of an axial isometry $\ga\in\Ga$. Then there exists a positive $\eps$ and a round 1-sphere $\si\subset\tits P(g)$ 
such that  $\dot B_\eps(g(+\infty))\cap\si$ 
does not contain a branch point. In particular, $\si\subset\tits X$ contains a non-empty open relatively compact subset.
\eprop

\proof
Since the diameter of $\tits X$ is $\pi$, the axis $g$ bounds a flat half-plane $H\subset X$.
If $\ga$ preserves $H$, then by Lemma~\ref{lem_per2fl}, the parallel set $P(g)$ contains a periodic 2-flat. 
In this case, the claim follows from Corollary~\ref{cor_finite_branch}.
If $H$ is not preserved, then $\si=\tau^+\cup\ga\tau^+$ is a round sphere in $\tits X$ where $\tau^+:=\geo H$.
Let $F$ be a 2-flat with $\geo F=\si$. Let $a>0$ be such that $g\subset N_a(F)$.
Now suppose that there is a sequence $(\xi_k^+)$ of pairwise distinct branch points in $\si$ with $\xi_k^+\to g(+\infty)$.
Denote by $\xi_k^-$ the antipode of $\xi_k^+$ in $\si$.
Since $\tits X$ has diameter $\pi$, we find for each $k\in\N$ a sequence of local geodesics $\al_{kl}\subset \tits X$ such that 
\begin{itemize}
	\item $\D\al_{kl}=\{\xi^-_{kl},\xi^+_{kl}\}$ and $\xi^\pm_{kl}\to\xi_k^\pm$ with respect to the Tits metric;
	\item the lengths of $\al_{kl}$  converges to $\pi$ as $l\to\infty$.
\end{itemize}
Hence by the \hyperref[lem_key_tech]{Half-Plane Lemma}, there exists a sequence of 2-flats $(\tilde F_{k})$ in $P(g)$, 
and a sequence of flat half-planes $(\tilde H_{k})$ with boundaries $\tilde h_k:=\D\tilde H_k\subset\tilde F_k$ and such that
the following properties hold.
\begin{enumerate}
	\item $g\subset N_a(\tilde F_k)$ for every $k\in\N$;
	\item $|\geo \tilde h_k,\geo g|=|\{\xi_k^-,\xi_k^+\},\geo g|$;
	\item $\geo\tilde H_k\cap\geo\tilde F_k=\geo\tilde h_k$.
\end{enumerate}
In particular, $\geo \tilde h_k\to \geo g$ with respect to the Tits metric.
Hence Lemma~\ref{lem_no_acc} implies that $\tilde h_k$ is parallel to $g$ for almost all $k$.
This is a contradiction since the points $\xi_k^+$ are pairwise distinct.
\qed

\subsection{The Tits boundary is a building}

\bdfn
For a CAT(0) space of rank 2, we define its {\em essential Tits boundary} $\etits X$ as the subset of $\tits X$ given by the union of all simple closed geodesics.
In particular, $\etits X$ is a 1-dimensional CAT(1) space which is geodesically complete. 
\edfn

\blem\label{lem_diameter_etits}
Let $X$ be a CAT(0) space of rank 2 which contains a 2-flat. If the diameter of $\tits X$ is equal to $\pi$, then so is the diameter of $\etits X$.
\elem

\proof
Since $X$ contains a 2-flat, the essential Tits boundary is non-empty. 
Let $\xi$ and $\hat\xi$ be points in $\etits X$. Then there is a geodesic $\al$ of length at most $\pi$ between them in $\tits X$.
We will show that this geodesic is contained in $\etits X$. Since $\etits X$ is geodesically complete, we may assume that $\al$ has length $\pi$
to begin with. Let $\eta$ be a point on $\al$ and let $\eps>0$ be such that $\eta\notin N(\{\xi,\hat\xi\})$.
Now we extend $\al$ beyond $\xi$ to a point $\xi'$ by a geodesic $\al'$ of length $\eps$ in $\etits X$.
Let $\beta$ be a geodesic in $\tits X$ from $\xi'$ to $\hat \xi$. 
Then $\hat\xi\notin\beta$ and therefore $\beta\cap(\al\setminus B_\eps(\xi))=\emptyset$.
The union $\al\cup\al'\cup\beta$ contains a simple closed geodesic $c$.
By definition $c$ has to lie in $\etits X$, and by construction $c$ contains $\al\setminus B_\eps(\xi)$ and therefore $\eta$. 
Thus, all of $\al$ is contained in $\etits X$ as required. 
\qed

\blem\label{lem_reg_1D}
Let $Z$ be a 1-dimensional CAT(1) space of diameter $\pi$.
Then the subset $O\subset Z$ of regular points is closed under taking antipodes.
\elem

\proof
Let $\xi\in O$ be a regular point and let $\hat\xi\in Z$ be an antipode.
Then there exists $s>0$ such that $\bar B_s(\xi)$ is isometric to a closed interval.
Denote by $\xi_k^\pm$ the two points at distance $\frac{s}{k}$ from $\xi$.
Note that we have the inclusion of geodesics $\xi_k^\pm\hat\xi\subset \xi_{k+1}^\pm\hat\xi$.
We conclude $|\xi,\xi^\pm|+|\xi^\pm,\hat\xi|=\pi$ and $\xi,\hat\xi$ lie in a round 1-sphere $\si$.
We claim that $B_s(\hat\xi)\subset\si$. Let $\hat\eta\in Z$ be  a point in  $B_s(\hat\xi)$.
We may assume that the geodesic $\hat\eta\xi$ contains the point $\xi_{1}^+$.
In particular, $|\hat\eta,\xi_{1}^+|\leq\pi-s$. Hence the geodesic triangle with vertices $\hat\eta, \xi_{1}^+$
and $\hat\xi$  has perimeter less than $2\pi$. Thus it is degenerated and we have $\hat\eta\in\xi_{1}^+\hat\xi$
as required. 
\qed

\blem\label{lem_etits_build}
Let $X$ be a locally compact CAT(0) space of rank 2 with a geometric group action. 
Suppose that the diameter of $\tits X$ is equal to $\pi$. Then $\etits X$ is a spherical join or a spherical building.
\elem

\proof
We may assume that $\etits X$ is not a round sphere. By \cite[Theorem~C]{Kleiner}, $X$ contains a 2-flat.
Thus, by  Lemma~\ref{lem_diameter_etits}, the diameter of $\etits X$ is equal to $\pi$. By \cite[Theorem~11]{Sw_cut}, the group $\Ga$
contains an axial element.
Hence, the open subset $O\subset\etits X$ of regular points is non-empty by Proposition~\ref{prop_regular_points}.
We infer from Lemma~\ref{lem_reg_1D} that $\etits X\setminus O$ is a proper closed subset of $\etits X$ which contains with every point all of its antipodes.
The claim follows from Theorem~\ref{thm_lytchak}.
\qed

\bcor\label{cor_rank2_build}
Let $X$ be a locally compact CAT(0) space of rank 2 with a geometric group action. 
Suppose that the diameter of $\tits X$ is equal to $\pi$. Then $\tits X$ is a spherical join or a spherical building.
\ecor

\proof
By Lemma~\ref{lem_etits_build}, it is enough to show that $\tits X=\etits X$. Suppose that  $\xi\in\tits X\setminus\etits X$
and $\xi'\in\etits X$ are points at distance less than $\pi$ such that the geodesic $\xi\xi'$ intersects $\etits X$ only in the point $\xi'$.
Let $\hat\xi'\in \etits X$ be an antipode of $\xi'$.
Let $\al$ be a geodesic from $\xi$ to $\hat\xi'$. Denote by $\eta$ the last point on  $\al\cap\xi\xi'$.
Since $|\xi',\eta|+|\eta,\hat\xi'|\geq\pi$ and $|\xi,\eta|+|\eta,\hat\xi'|\leq\pi$, we see $|\xi',\eta|\geq |\xi,\eta|$.
Hence the union of $\xi'\xi\cup\al$ with a geodesic $\beta$ in $\etits X$ joining $\xi'$ to $\hat\xi'$ contains a simple closed geodesic which covers
at least half to the segment $\xi'\xi$. By the definition of $\etits X$, and the choice of $\xi'$,
we arrive at a contradiction.
\qed

\proof[Proof of Theorem~\ref{thm_mainB}]
By Corollary~\ref{cor_rank2_build}, the Tits boundary $\tits X$ is a spherical join or a building.
Since $X$ is geodesically complete and locally compact, the claim follows from Theorem~\ref{thm_leeb} and \cite[Proposition~2.3.7]{KleinerLeeb}.
\qed

\appendix
\section{The relative asymptotic Plateau problem after Kleiner--Lang}\label{sec_app}

The primary goal of this section is to provide a proof of the \hyperref[lem_key_tech]{Half-Plane Lemma}.
Our proof requires a solution of an asymptotic Plateau problem relative to a flat.
Since this has not been done in the literature, we need some preparation.
The (non-relative) asymptotic Plateau problem has recently been solved  by Kleiner--Lang~\cite{KL_higher}
and we use their work extensively throughout this section. The idea to use solutions to a relative Plateau problem originated from \cite{HKS_I}.

In Section~\ref{subsec_notation}, we agree on notation and proof  monotonicity and volume rigidity of minimizing currents in CAT(0) spaces,
a result which is folklore.
In Section~\ref{subsec_relative}, we solve the asymptotic Plateau problem relative to a 2-flat.
We provide basic properties of relative minimizers which are required for our proof of the \hyperref[lem_key_tech]{Half-Plane Lemma}.
At last, in Section~\ref{subsec_rel_asym}, we present the proof of the \hyperref[lem_key_tech]{Half-Plane Lemma}.
We have refrained from any kind of general treatment, in particular we only discuss spaces of rank 2. 
This is all we need here and a broader discussion would get out of hand. 

The general reference for this section is \cite{KL_higher}.
See also \cite{AK, La_currents} for a thorough background on geometric measure theory in metric spaces.

\subsection{Currents in CAT(0) spaces}\label{subsec_notation}

Let $X$ be a locally compact metric space.
For every integer $n \ge 0$, let $\cD^n(X)$ denote the set of all 
$(n+1)$-tuples $(\pi_0,\ldots,\pi_n)$ of real valued functions on $X$ such 
that $\pi_0$ is Lipschitz with compact support $\spt(\pi_0)$ and 
$\pi_1,\dots,\pi_n$ are locally Lipschitz.
An {\em $n$-current} $S$ in $X$ is a function 
$S:\cD^n(X) \to \R$ satisfying the following three conditions:
\begin{enumerate}
	\item(multilinearity)
$S$ is $(n+1)$-linear;
\item(continuity) 
$S(\pi_{0,k},\ldots,\pi_{n,k}) \to S(\pi_0,\ldots,\pi_n)$
whenever $\pi_{i,k} \to \pi_i$ pointwise on $X$, $\bigcup_k\spt(\pi_{0,k})$ is bounded 
and the $\pi_{i,k}$ are uniformly locally Lipschitz continuous; 
\item(locality)
$S(\pi_0,\ldots,\pi_n) = 0$ whenever one of the functions
$\pi_1,\ldots,\pi_n$ is constant on a neighborhood of $\spt(\pi_0)$.
\end{enumerate}

We write $\cD_n(X)$ for the vector space of $n$-currents 
in $X$. For every $S \in \cD_n(X)$ we denote by  $\spt(S)$ the {\em support} of $S$ and by  $\D S \in \cD_{n-1}(X)$ its {\em boundary}.
As usual, $\|S\|$ will denote the associated regular Borel measure and $\M(S)=\|S\|(X)$ its {\em mass}. 
The {\em flat norm} of $S$ is denoted $\cF(S)$.
Moreover, for a Borel set $A\subset X$ we denote by $S\on A \in \cD_n(X)$ the {\em restriction} of $S$ to $A$.
For a proper Lipschitz map $f: X \to Y$ into another locally compact
metric space $Y$, we have the {\em push-forward} $f_\#S \in \cD_n(Y)$.
We denote by $\bI_{n,loc}(X)$ the abelian group of {\em locally integral $n$-currents},
and write $\bI_{n,c}(X)$ for the subgroup of {\em integral $n$-currents with compact support}.
Further, $\bZ_{n,loc}(X)\subset\bI_{n,loc}(X)$ and $\bZ_{n,c}(X)\subset\bI_{n,c}(X)$
will denote the subgroup of {\em cycles}, i.e. integral currents with boundary  zero.

Recall that an integral current $S\in \bI_{n,c}(X)$ is  concentrated on a countable $\Ha^n$-rectifiable set $\chi_S$, its {\em characteristic set} 
\cite[Theorem~4.6]{AK}. The characteristic set $\chi_S$ is unique in the sense that any Borel set $E$ with $\|S\|=\|S\|\on E$
contains $\chi_S$ up to $\Ha^n$-negligible sets. Moreover, there exists a $\Ha^n$-integrable {\em multiplicity function} $\theta:\chi_S\to\N$
such that $\|S\|=\theta\cdot\Ha^n\on\chi_S$ \cite[Theorem~9.5]{AK}. Note that by \cite[Lemma~9.2]{AK} the area factor $\la$ appearing in the general formula in 
\cite[Theorem~9.5]{AK} is equal to $1$ in our case, since CAT(0) spaces have Euclidean tangent cones, cf.~\cite[Section~11]{LW_plateau}.

On  Euclidean space $\R^m$ we denote the {\em flat $n$-chains with compact support} 
by $\cF_{n,c}(\R^m)$, see \cite{Fl_flat, W_defo, W_rec}.
Recall that for a compact subset $K\subset\R^m$ the space of flat chains with compact support in $K$ is the flat-closure
of $\cP_{n}(K)$, {\em the polyhedral $n$-chains with support in $K$}.
We will need a localized version. Define
 a {\em local flat $n$-chain} as an $n$-current $S$
such that for every $x\in\R^m$ there exists $S'\in\cF_{n,c}(\R^m)$ with $x\notin\spt(S-S')$.
The space of  {\em locally flat $n$-chains} is denoted by $\cF_{n,loc}(\R^m)$.

Let  $X$ be a CAT(0) space and $S\in \bZ_{n,c}(X)$. Then for every point $p\in X$ we have the {\em cone}
from $p$ over $S$, denoted by $C_p S\in\bZ_{n+1,c}(X)$, cf. \cite[Section~2.7]{KL_higher}. If $\spt(S)\subset \bar B_R(p)$ the following
{\em Euclidean cone inequality} holds.
\[\M(C_p S)\leq \frac{R}{n}\cdot\M(S).\]  

The fundamental class of $\R^n$ is denoted by $\bb{\R^n}\in\bZ_{n,loc}(\R^n)$. Every proper Lipschitz map $\varphi:\R^n\to X$
induces a natural local cycle $\varphi_\#(\bb{\R^n})\in\bZ_{n,loc}(X)$. In this way, a flat $F$ in $X$ becomes a local cycle which we will still denote by $F$
and call a {\em multiplicity 1 flat}.

Following \cite{KL_higher}, for a current
$S \in \bI_{n,\loc}(X)$, a point $p \in X$ and $r > 0$, we define the {\em ($r$-)density at $p$} by 
\[
\G_{p,r}(S) := \frac{1}{r^n} \|S\|(B_r(p))\quad \G_{p}(S) :=\liminf\limits_{r\to 0} \frac{1}{r^n} \|S\|(B_r(p)).
\]  
 Furthermore, for any $p \in X$, we define the  {\em density at infinity} by  
\[
\Gi(S) := \limsup\limits_{r\to\infty} \G_{p,r}(S).
\]
Similarly, we define the {\em filling density at $p$} by 
\[
\F_{p,r}(S) := \frac{1}{r^{n+1}} 
\inf\{\M(V) : V \in \bI_{n+1,c}(X),\,\spt(S-\D V) \cap B_r(p) =\emptyset\}.
\]  
Furthermore, for any $p \in X$, we define the  {\em filling density at infinity} by 
\[
\Fi(S) := \limsup\limits_{r\to\infty} \F_{p,r}(S).
\]

The following monotonicity property of minimizers is certainly well-known. We include
a version for CAT(0) spaces, since we were unable to find a reference. 

\bprop[monotonicity]\label{prop_mon}
Let $X$ be a locally compact CAT(0) space and $S \in \bI_{n,c}(X)$ an area minimizer.
Then for every point $p\in\spt(S)$  the $r$-density
\[\om_n\leq \G_{p,r}(S)\]
is a non-decreasing function of $r$ as long as $r\leq|p,\spt(\D S)|$.
Moreover, if $\G_{p,r_0}(S)=\om_n$ holds for some $r_0\leq|p,\spt(\D S)|$, then $\spt(S)\cap \bar B_{r_0}(p)$ is 
isometric to a Euclidean ball of radius $r_0$, and $S$ has constant multiplicity 1.
\eprop

\proof
Monotonicity follows from the standard cone comparison since CAT(0) spaces satisfy a Euclidean coning inequality.
More precisely, for almost all $r\leq |p,\spt(\D S)|$ holds 
\[\M(S\on B_r(p))\leq \frac{r}{n}\cdot\M(\D(S\on B_r(p))).\]
By the coarea inequality we have 
\[\frac{d}{dr}\M(S\on B_r(p))\geq \M(\D(S\on B_r(p)))\] 
for almost all $r$, and monotonicity follows by integration.
Recall that $S$ is concentrated on the countably $\Ha^n$-rectifiable characteristic set $\chi_S$ which is given by $\chi_S:=\{x\in X|\ \G_x(S)>0\}$.
The isoperimetric inequality yields some positive lower density bound at all points $p\in\spt(S)$ \cite[Lemma~3.3]{KL_higher}.
This implies $\chi_S=\spt(S)$, so $\chi_S$ is closed.
It follows from \cite[Theorem~5.4]{AK_rec} that $\G_p(S)=\theta\cdot\om_n$ for almost all points $p\in\chi_S$ where $\theta:\chi_S\to \N$ denotes the multiplicity function of $S$.
The bound  $\G_p(S)\geq\om_n$ extends to the support $\spt(S)$ by upper semi-continuity of $\G_p(S)$.

Now suppose $\G_{p,r_0}(S)=\om_n$ holds for some $r_0\leq|p,\spt(\D S)|$.
Then the above inequalities become equalities and we have
\[\M(S\on B_r(p))\equiv \om_n\cdot r^n\quad \text{and}\quad \M(\D(S\on B_r(p)))\equiv n\cdot\om_n\cdot r^{n-1}.\]
Set $g:=(d_p)|_{\chi_S}$. 
By the general coarea formula \cite[Theorem~9.4]{AK_rec}, 
the coarea factor $C_1(d^S g_x)$ has to be $1$ almost everywhere on $\chi_S\cap \bar B_{r_0}(p)$, cf.~\cite[Section~9]{AK_rec}. 
Let $x\in\chi_S$ be such that $\chi_S$ has a tangent space $T_x\chi_S\subset T_x X$\footnote{
While in the general setting of \cite{AK_rec} one has to work with approximate tangent spaces, in the CAT(0) setting the 
Lytchak-Rademacher theorem \cite[Theorem~1.6]{Ly_diff} allows us
to realize the approximate tangent space of a countably $\Ha^n$-rectifiable subset in $X$ as a linear subspace of the tangent space to $X$.} 
and $g$ has a tangential differential $d^S_x$ at $x$ \cite[Theorem~8.1]{AK_rec}. Let $\{v_1,\ldots,v_n\}$ be an orthonormal basis 
for $T_x \chi_S$ such that $v_1$ is the nearest point to $\log_x(p)\in T_x X$. Then, by the first variation formula, 
$d^S g_x=(-\cos(\al_1),0,\ldots,0)$ where $\al_1=\angle_x(\log_x(p),v_1)$. Thus, $C_1(d^S g_x)=|\cos(\al_1)|$. 
Hence,
the tangent space $T_x\chi_S$ contains the direction $\log_x(p)\in T_x X$ at $\Ha^n$-almost all points $x\in\chi_S\cap \bar B_{r_0}(p)$. For $r>0$
denote by $\pi_r:X\setminus B_r(p)\to S_r(p)$ the nearest point projection.
Then, by the general area formula \cite[Theorem~8.2]{AK_rec}, for every $r<r_0$ the intersection 
$C_p(\chi_S\cap \bar B_{r_0}(p))\cap S_r(p)=\pi_r(\chi_S\cap \bar B_{r_0}(p))$ is $\Ha^n$-negligible.
Hence the cone $C_p(\chi_S\cap\bar B_{r_0}(p))$ is $\Ha^{n+1}$-negligible and therefore $C_p(S\on \bar B_{r_0}(p))=0\in\bI_{n+1,c}(X)$.
Since $\D C_p(S\on \bar B_{r_0}(p))=S\on \bar B_{r_0}(p)-C_p(\D (S\on \bar B_{r_0}(p)))$, we conclude that $S$ is conical with respect to $p$, $S\on \bar B_{r_0}(p)=C_p(\D (S\on \bar B_{r_0}(p)))$. In particular, if $c_\la:\bar B_{r_0}(p)\to \bar B_{\la r_0}(p)$ the $\la$-Lipschitz map which 
contracts radial geodesics relative $p$ by the factor $\la\in(0,1)$, then we have $(c_\la)_\# S\on \bar B_{r}(p)=S\on \bar B_{\la r}(p)$ for $r\leq r_0$.
Moreover, since the characteristic set is closed, $\chi_S\cap \bar B_{r_0}(p)$ agrees with $C_p(\chi_S\cap \bar B_{r_0}(p))$.
By monotonicity and the mass control, the multiplicity function $\theta$ has to be equal to $1$ in a neighborhood of $p$ in $\chi_S$.
Since $S$ is conical, $\theta\equiv 1$ on all of $\chi_S\cap \bar B_{r_0}(p)$. This implies for all $r<r_0$, 
\[\Ha^n(\chi_S\cap B_r(p))=\om_n\cdot r^n\quad \text{and}\quad \Ha^{n-1}(\chi_S\cap S_r(p))= n\cdot\om_n\cdot r^{n-1}.\]
Thus $c_\la$ scales $\Ha^n$-measure on $\chi_S\cap B_{r_0}(p)$ by the factor $\la^n$.
We claim that for every $x\in\spt(S)$ with $a=|p,x|$ and $r\leq r_0-a$ we have $\|S\|(\bar B_{r}(x))=\om_n\cdot r^n$.
Indeed, using the lower density bound at $x$, monotonicity and the effect of $c_\la$ on $\Ha^n$, we obtain
\begin{align*}
\om_n\cdot r^n&\leq \|S\|(B_r(x))=\la^{-n}\cdot \|S\|(c_\la(B_r(x)))\\
&\leq\la^{-n}\cdot \|S\|(B_{\la r}(c_\la(x)))\\
&\leq \la^{-n}\cdot\left(\frac{\la r}{r_0-\la a}\right)^n\cdot\om_n\cdot r_0^n
\stackrel{\la\to 0}{\longrightarrow} \om_n\cdot r^n.
\end{align*}
We infer $\chi_S\cap\bar B_r(x)=C_x(\chi_S\cap\bar B_r(x))$. Therefore $\chi_S\cap\bar B_{r_0}(p)$ is conical with respect to all of its points.
In particular, it is a closed convex subset of $X$. From the mass bound, it follows that  $\chi_S\cap\bar B_{r_0}(p)$ is isometric to an Euclidean $r_0$-ball.
\qed

From the rigidity in Proposition~\ref{prop_mon}, we obtain:

\bcor\label{cor_mon}
Let $X$ be a locally compact CAT(0) space and $S \in \bI_{n,loc}(X)$ an area minimizer.
If $\Gi(S)=\om_n$, then $\spt(S)$ is an $n$-flat in $X$ and $S$ has constant multiplicity 1.
\ecor

\subsection{The asymptotic Plateau problem relative to a 2-flat}\label{subsec_relative}

We fix the following setup for this section. Let $X$ be a locally compact CAT(0) space of rank 2 which contains a round sphere 
$\si\subset\tits X$ in its Tits boundary. We fix a 2-flat $F\subset X$ with $\tits F=\si$ and a base point $o\in F$. 
Moreover, $\al^+\subset\tits X$ denotes an embedded local geodesic with $\al^+\cap\si=\D\al^+$.
Our goal is to solve an asymptotic Plateau problem relative $F$. Informally speaking, we search for a minimal surface $S\subset X$
whose boundary lies in $F$ and which is asymptotic to $\al^+$. If the length of $\al^+$ is close to $\pi$, then we expect $S$
to behave almost like a flat half-plane orthogonal to $F$.
In order to achieve this, we employ doubling to reduce the problem to the ordinary asymptotic Plateau problem treated in \cite{KL_higher}. 
So we consider the double of $X$ along $F$:
\[\hat X=X^-\cup_F X^+\]
where $X^\pm$ are two isometric copies of $X$.  Note that $\tits \hat X=\tits X^-\cup_\si \tits X^+$.
In particular, $\hat X$ has rank 2. Let us denote by $\iota$ the natural isometry of $\hat X$ which interchanges $X^-$
and $X^+$. For subsets $A^+\subset X^+$ we set $A^-:=\iota (A^+)\subset X^-$. We call a current $S\in\bI_{2,loc}(\hat X)$
{\em symmetric}, if 
\[S=S^+-\iota_\#(S^+)\]
where $S^+:=S\on X^+$. Finally, we also double our boundary data, and define $\al:=\al^+\cup\iota(\al^+)\subset\tits X$.

\bdfn
For $S\in \bI_{n,loc}(\hat X)$  set $S^+:=S\on X^+$. Then $S^+$ is a {\em relative minimizer} if 
\[\M(S^+\on B^+)\leq \M(W^+)\]
whenever $B^+\subset X^+$ is a Borel set such that $S\on (B^+\cup B^-)\in\bI_{n,c}(\hat X)$ and $W^+\in\bI_{n,c}(X^+)$
satisfies $\spt(\D(S^+\on B^+-W^+))\subset F$.
\edfn

\blem\label{lem_flat_chain}
Let $S^+\in\cD_{2,loc}(X^+)$ be a current of locally finite mass and with $\spt(\D S^+)\subset F$.
Suppose that  $\|S^+\|(F)=0$ and that there is a sequence $t_k\to 0$, such that $S^+\on(X\setminus N_{t_k}(F))\in\bI_{2,loc}(X^+)$.
Then $\D S^+$ is a local flat chain in $F$, $\D S^+\in\cF_{1,loc}(F)$. Moreover, the double $S:=S^+-\iota_\#(S^+)$ is a locally integral cycle,
$S\in\bZ_{2,loc}(\hat X)$.
\elem

\proof
For a  point $p\in X^+$ we denote by $d_p$ the distance function to $p$.
To see that $\D S^+$ is a local flat chain in $F$, it is enough to prove the following claim.
For any $p\in\spt(\D S^+)$ there exists $r_0>0$ such that $\<S^+\on(X\setminus N_{t_k}(F)),d_p,r_0\>\in\bI_{1,c}(X^+)$ for every 
$k\in\N$, and $S^+\on\bar B_{r_0}(p)$ is a flat limit of integral chains $T_k\in\bI_{2,c}(X^+)$.

Indeed, if the claim holds, then, since the nearest point projection $\pi_F:X^+\to F$ is Lipschitz, $(\pi_F)_\#(S^+\on\bar B_{r_0}(p))\in\cF_{2,c}(F)$
and therefore $\D(\pi_F)_\#(S^+\on\bar B_{r_0}(p))=(\pi_F)_\#\D(S^+\on\bar B_{r_0}(p))\in\cF_{1,c}(F)$.
Then, either $p\notin\spt(\D S^+ -(\pi_F)_\#\D(S^+\on\bar B_{r_0}(p)))$ or, there exists a large $k\in\N$ such that 
$p\notin\spt(\D S^+ -(\pi_F)_\#\<S^+\on(X\setminus N_{t_k}(F)),d_p,r_0\>$. In any case, we see that $\D S^+$ is locally flat near $p$.

To prove the claim, we set $T_k:=S^+\on(X\setminus N_{t_k}(F))\in\bI_{2,loc}(X^+)$. By continuity of $S^+$ as an extended functional 
\cite[Theorem~4.4]{La_currents}, we have $T_k\to S^+\on(X^+\setminus F)$ weakly.
Since $\|S^+\|(F)=0$, we have $S^+\on(X^+\setminus F)=S^+$ \cite[Lemma~4.7]{La_currents}.
For $r>0$ we infer from \cite[Theorem~6.2 (3)]{La_currents}
\[\int_0^\infty \|\<T_k,d_p,s\>\|(B_r(p))\ ds\leq\|T_k\|(B_r(p))\leq\|S^+\|(B_r(p)).\]
By monotone convergence, we conclude
\[\int_0^\infty\lim\limits_{k\to\infty} \|\<T_k,d_p,s\>\|(B_r(p))\ ds\leq\|S^+\|(B_r(p)).\]
Now choose $r_0\in(0,r)$ such that $\<T_k,d_p,r_0\>\in\bI_{1,c}(X^+)$ for every $k\in\N$ and $\lim\limits_{k\to\infty} \|\<T_k,d_p,r_0\>\|<\infty$.
Set $T_k':=T_k\on\bar B_{r_0}(p)$. Because $\|S^+\|(F)=0$, we see $\lim\limits_{t\to 0}\|S^+\|(\bar B_{r_0}(p)\cap N_t(F))=0$. Then
for $k<l$ we have 
\[\cF(T'_k,T'_l)\leq (\|S^+\|+\lim\limits_{k\to\infty} \|\<T_k,d_p,r_0\>\|)(\bar B_{r_0}(p)\cap \bar N_{t_k}(F))\stackrel{k,l\to \infty}{\rightarrow}0.\]
Hence $(T'_k)$ is a Cauchy sequence with respect to flat distance.
Because flat convergence implies weak convergence, the limit coincides with $S^+\on\bar B_{r_0}(p)$. This confirms the claim and therefore $\D S^+\in\cF_{1,loc}(F)$.

By the deformation theorem \cite{FF_currents,W_defo}, we can now find a locally polyhedral chain $P\in\cP_{1,loc}(F)$ and a locally flat 
chain $V\in\cF_{2,loc}(F)$ such that $\D V=P-\D S^+$. Recall that since $V$ is top-dimensional in $F$, it is canonically identified with a locally integrable
function with integer values. In particular, $S^+ +V$ is locally integral. Hence $S=(S^+ +V)-\iota_\#(S^+ +V)$ is a locally integral cycle as claimed.
\qed

\bcor\label{cor_flat_chain}
Let $S\in\bZ_{2,loc}(\hat X)$ be symmetric, $S=S\on X^+-\iota_\#(S\on X^+)$.
Then the boundary of the restriction $S^+:=S\on X^+$ is a local flat chain in $F$, $\D S^+\in\cF_{1,loc}(F)$.
\ecor

\proof
Note that since $S$ is symmetric, we must have $\|S\|(F)=0$. Hence the claim follows from slicing and Lemma~\ref{lem_flat_chain}.
\qed

\blem\label{lem_minimizer}
There exists a minimizing $S\in\bZ_{2,loc}(\hat X)$  which is symmetric, and such that 
\[\Fi(S-C_o\al)=0\quad \text{ and }\quad \Gi(S)=\frac{\length(\al)}{2}.\] 
Moreover, the restriction $S^+:=S\on X^+$
is a relative minimizer.
\elem

\proof
The existence of a minimizer follows from \cite[Theorem~5.6]{KL_higher}. 
We only need to argue that we can choose it to be symmetric. As in \cite[Theorem~5.6]{KL_higher}, we choose a sequence of radii $r_k\to \infty$ and solve the Plateau Problem
for $C_o(\al)\on B_{r_k}(o)$ to obtain minimizers $S_k\in\bI_{2,c}(\hat X)$ with $\D S_k=\D (C_o(\al)\on B_{r_k}(o))$.
Note that since $C_o(\al)\on B_{r_k}(o)$ is symmetric, we may assume that $S_k$ is symmetric as well. Indeed, by Lemma~\ref{lem_flat_chain},  
$\tilde S_k:=(S_k\on X^+)-\iota_\#(S_k\on X^+)$ lies in $\bI_{2,c}(\hat X)$ and $\M(\tilde S_k)\leq 2\cdot\M(S_k\on X^+)\leq\M(S_k)$. 
Now the proof in  \cite[Theorem~5.6]{KL_higher} applies to produce a symmetric minimizer $S\in\bZ_{2,loc}(\hat X)$.
By \cite[Theorem~8.3]{KL_higher}, $S$ is $\F$-asymptotic to $C_o\al$ and has the required volume growth. That $S^+$ is a 
relative minimizer follows from symmetry of $S$.
\qed

\medskip

In Section~\ref{subsec_no_acc} we used flat half-planes orthogonal to $F$ to rule out accumulation of branch points at infinity.
Just like the relative minimizer $S^+$ serves as a substitute for a flat half-plane, the following lower bound on filling density serves as a substitute
for orthogonality to $F$.

From now on let $S\in\bZ_{2,loc}(\hat X)$ be a symmetric minimizer as in Lemma~\ref{lem_minimizer}.

\blem\label{lem_lower_dist_bound}
For every $p\in\spt(S)\cap F$ and $r>0$ holds
\[\F_{p,r}(S-F)\geq\frac{\pi}{3}.\]
\elem

\proof
Let $V\in\bI_{3,c}(\hat X)$ be such that $\spt(S-F-\D V)\cap B_r(p)=\emptyset$.
For almost all $s\in[0,r]$ we have $\<S,d_p,s\>-\<F,d_p,s\>=\D\<V,d_p,s\>$ and all occurring slices are integral. Since $S$ is a symmetric minimizer, we must have
\[\|\<V^\pm,d_p,s\>\|(B_s(p))\geq\|S^\pm(B_s(p))\|.\]
We obtain $\|S^\pm(B_s(p))\|\geq\frac{\pi}{2}\cdot s^2$ from Proposition~\ref{prop_mon}. Since 
\[\M(V)\geq 2\cdot\M(V^+)\geq 2\cdot\int_0^r \M(\<V^+,d_p,s\>)\ ds,\]
the claim follows by integration.
\qed

\blem\label{lem_rel_visibility}
For every $\eps>0$ there exists $r_\eps>0$
such that for all $r\geq r_\eps$ holds
\[\spt(\D S^+)\cap B_r(o)\subset N_{\eps r}(C_o(\D\al^+)).\]
\elem

\proof
Suppose for contradiction that there exists $\eps_0>0$ and a sequence $x_k\in\spt(\D S^+)$
with $|o,x_k|\to\infty$ but $|x_k,C_o(\D\al^+)|\geq\eps_0\cdot|o,x_k|$. 
Since $\spt(\D S^+)\subset\spt(S^+)\subset\spt(S)$, we infer from \cite[Theorem~8.1]{KL_higher}, that there exist points $y_k\in C_o(\al)$ with $\frac{|x_k,y_k|}{|x_k,o|}\to 0$.
After passing to a subsequence, we find an element $\eta\in\al$ such that $\frac{|y_k,o\eta|}{|y_k,o|}\to 0$.
Hence $\frac{|x_k,o\eta|}{|x_k,o|}\to 0$. But since $x_k\in F$, we have $\eta\in\D\al^+$. Contradiction.
\qed

\blem\label{lem_intersection}
Let $\zeta$ and $\hat\zeta$ be antipodes on $\si$
which separate $\D\al^+$. Then $\spt(\D S^+)$ intersects every complete geodesic $g_0\subset F$ with $\geo g_0=\{\zeta,\hat\zeta\}$.
\elem

\begin{center}
\includegraphics[scale=0.5,trim={-3cm 0cm 0cm 0cm},clip]{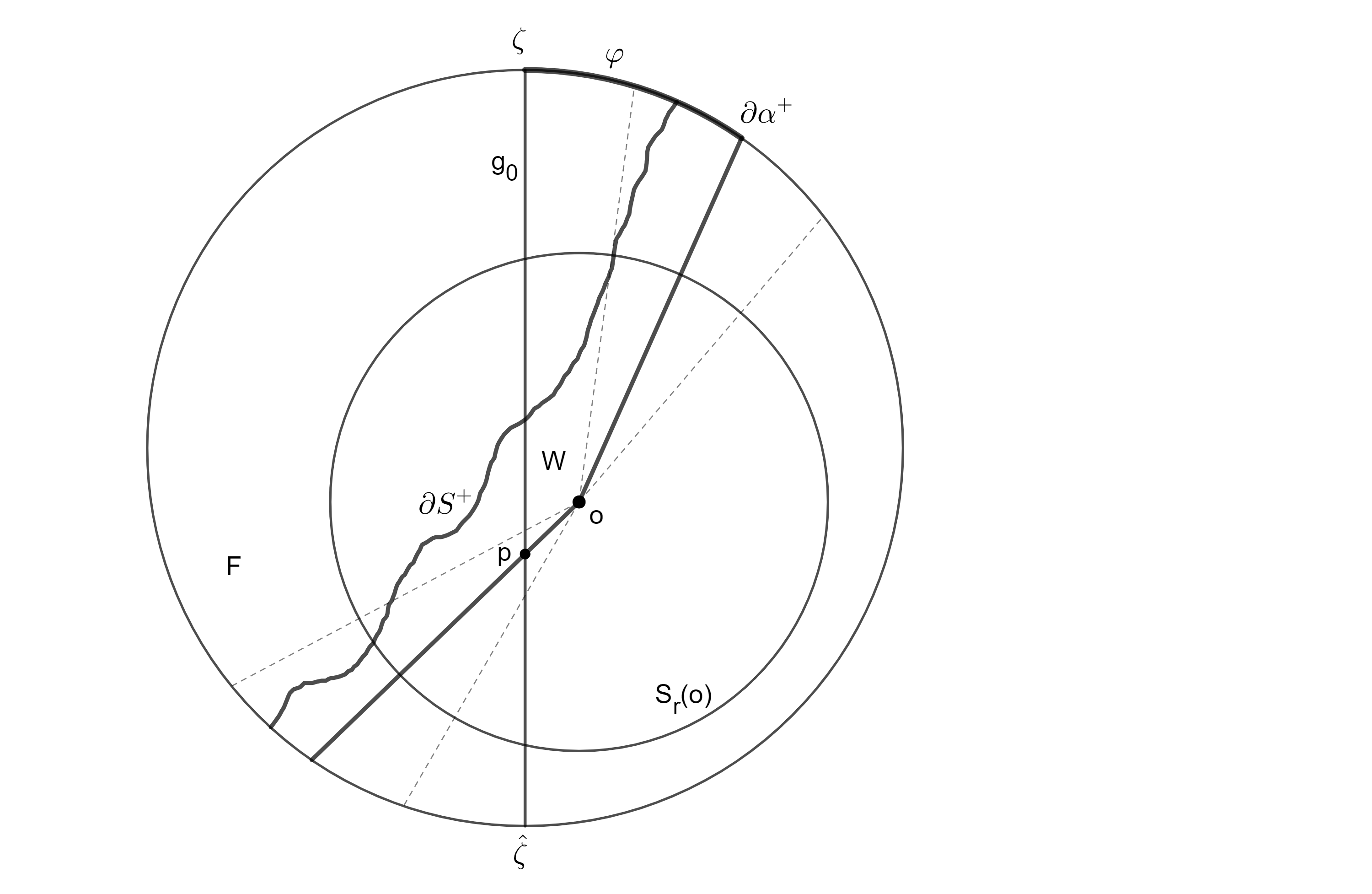}
\end{center}

\proof
By Lemma~\ref{lem_flat_chain}, $\D S^+$ is a local flat chain in $F$.
Note that $\Fi(S-C_o(\al))=0$ implies $\Fi(\D S^+-C_o(\D\al))<\infty$. 
Indeed, for $\eps_0>0$ choose $r>0$ and $V\in\bI_{3,c}(\hat X)$ such that 
$\spt(S-C_o\al-\D V)\cap B_r(o)=\emptyset$ and $\M(V)\leq\eps_0\cdot r^3$.
Denote by $d_F$ the distance function to $F$ in $X^+$ and by $\pi_F:\hat X\to F$ the nearest point projection.
By the coarea inequality, we find $s\in(0,\frac{r}{2})$ such that $\M(\<V,d_F,s\>)\leq\frac{2}{r}\cdot\M(V)\leq 2\eps_0\cdot r^2$
and $\spt(\<S-C_o\al,d_F,s\>-\D\<V,d_F,s\>)\cap B_r(o)=\emptyset$, and all occurring slices are integral.
Set $T:=(S^+-C_o\al^+)\on\bar N_s(F)-\D\<V,d_F,s\>\in\bI_{2,c}(X^+)$. Then 
$\spt(\D S^+-C_o\D\al^+-\D(\pi_F)_\# (T))\cap B_{\frac{r}{2}}(o)=\emptyset$ and $\M((\pi_F)_\# (T))\lesssim r^2$.

In particular, $\D S^+$ and $C_o(\D\al^+)$ are homologous as local flat chains.
Denote by $W\in\cF_{2,loc}(F)$ the canonical filling of $\D S^+ -C_o(\D\al^+)$. By Lemma~\ref{lem_rel_visibility},
the support of $W$ lies sublinearly close to $C_o(\D\al^+)$.
Let $\varphi$ be a positive lower bound for the distances between the sets $\{\zeta,\hat\zeta\}$ and $\D\al^+$ in $\tits F$.
Choose $\eps\leq\frac{\varphi}{2}$ and then choose $r_\eps>0$ according to Lemma~\ref{lem_rel_visibility}.
By assumption, $C_o(\D\al^+)$ intersects $g_0$ transverse in exactly one point, say $p\in F$.
Now we can find $r\geq \max\{r_\eps,2|o,p|\}$ such that all points in $\spt(\<W,d_o,r\>)$ have distance larger than $\eps\cdot r$ from the geodesic $g_0$.
 Since $\D(W\on B_r(o))=\<W,d_o,r\>+(\D S^+-C_o(\al^+))\on B_r(o)$ is a flat 1-cycle with compact support,
$\spt(\D S^+)$ has to intersect $g_0$ in $B_r(o)$. 
\qed

\subsection{Producing flat half-planes as relative minimizers}\label{subsec_rel_asym}

Now we have all the necessary ingredients to prove:

\begin{namedlemma}[Half-Plane]
Let $X$ be a locally compact CAT(0) space of rank 2.  Let $\si\subset\tits X$ be a round sphere.
Further, let $g$ be an axis of an axial isometry $\ga$ and assume $\geo g\subset\si$.
Let $\xi^\pm\in\si$ be antipodes disjoint from $\geo g$.
Suppose that there is a sequence of local geodesics $\al^+_k\subset \tits X$ such that 
\begin{itemize}
	\item $\D\al^+_k=\{\xi^-_k,\xi^+_k\}$ and $\xi^\pm_k\to\xi^\pm$ with respect to the Tits metric;
	\item the lengths of $\al_k^+$  converges to $\pi$ as $k\to\infty$.
\end{itemize}
Then there exists a 2-flat $\tilde F\subset P(g)$, and a flat half-plane $\tilde H$ with boundary $\tilde h:=\D \tilde H\subset\tilde F$ such that the following properties hold.
\begin{enumerate}
	\item $g\subset N_a(\tilde F)$ if $g\subset N_a(F)$ for $a>0$;
	\item $|\geo \tilde h,\geo g|=|\{\xi^-,\xi^+\},\geo g|$;
	\item $\geo\tilde H\cap\geo\tilde F=\geo\tilde h$.
\end{enumerate}
\end{namedlemma}

\begin{center}
\includegraphics[scale=0.5,trim={-1cm 0cm 0cm 0cm},clip]{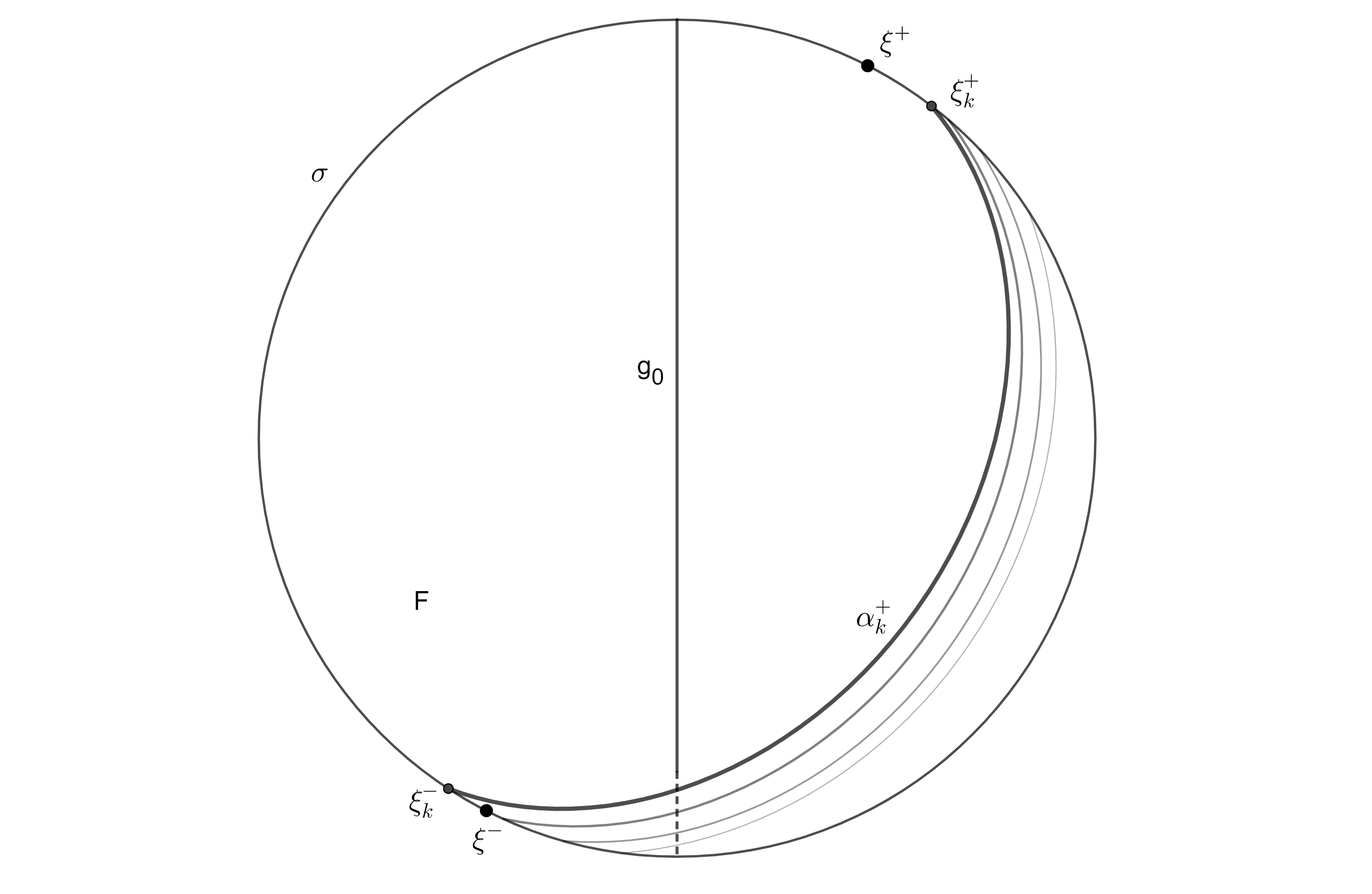}
\end{center}

\proof
Choose a 2-flat $F$ with $\geo F=\si$ and a base point $o\in F$. 
As before, we consider the double of $X$ along $F$, $\hat X=X^-\cup_F X^+$
and put $\al_k:=\al_k^+\cup \iota(\al_k^+)$.
By Lemma~\ref{lem_minimizer}, we find a symmetric minimizer $S_k\in\bZ_{2,loc}(\hat X)$ with $\Fi(S_k-C_o(\al_k))=0$
and $\Gi(S_k)=\frac{\length(\al_k)}{2}$. Let $g_0$ be a complete geodesic in $F$ parallel to $g$. By Lemma~\ref{lem_intersection},
there exists  a point $p_k\in\spt(S_k)\cap g_0$. 

Now let us choose numbers $m_k\in\Z$ such that the points $\ga^{m_k}(p_k)$ stay in a fixed compact set. Put $\ga_k:=\ga^{m_k}$.
After passing to subsequences, we have convergences with respect to pointed Hausdorff topology 
\[\ga_k(C_{p_k}(\al^+_k),p_k)\to(\tilde H,p_\infty)\quad\text{ and }\quad\ga_k (F,p_k)\to(\tilde F,p_\infty)\]
where $\tilde H$ is a flat half-plane, $\tilde F$ is a 2-flat, and $\tilde h:=\D\tilde H\subset\tilde F$.
We will now show that $\tilde F$ and $\tilde H$ satisfy the three required properties.
\medskip

\noindent (1) Since $\ga$ preserves $g$, we have $g\subset N_a(\tilde F)$ for $a>0$ whenever $g\subset N_a(F)$.
\medskip
	
\noindent (2) Let $\geo\tilde h=\{\tilde\xi^-,\tilde\xi^+\}$. 
Since $\xi^\pm_k\to\xi^\pm$ with respect to the Tits metric, and since $\ga$ preserves $g$, we have
$|g(\pm\infty),\ga_k\xi_k^+|=|g(\pm\infty),\xi_k^+|$ and $|g(\pm\infty),\ga_k\xi_k^-|=|g(\pm\infty),\xi_k^-|$.
Hence the statement follows from lower semi-continuity of the Tits metric with respect to the cone topology
\cite[Lemma~2.3.1]{KleinerLeeb}.
\medskip

\noindent (3) Note that if this fails, then  we must have $\tilde H\subset \tilde F$ since $\D\tilde H\subset\tilde F$.

Denote by $\tilde F^\pm\subset\tilde F$ the two half-planes determined by $\D\tilde H$.
Further, denote by $\si_k^\pm\subset\si$ the two closed arcs determined by $\{\xi_k^-,\xi_k^+\}$.
Then, as currents, we can write $F=F_k^+-F_k^-$ where $F_k^\pm=C_{p_k}(\si_k^\pm)\in\bI_{2,loc}(\hat X)$.
From \cite[Proposition~4.5]{KL_higher}, we conclude that for every $\eps>0$
there exists $a_\eps'>0$ such that $\F_{p_k,r}(S_k-C_{p_k}(\al_k))\leq\frac{\eps}{3}$ for all $r\geq a_\eps'$.
On the other hand, Lemma~\ref{lem_lower_dist_bound} implies $\F_{p_k,r}(S_k-F)\geq\frac{\pi}{3}$ for all $r>0$.
Hence, $\F_{p_k,r}(F-C_{p_k}(\al_k))\geq\frac{\pi-\eps}{3}$ and thus, by symmetry of  $C_{p_k}(\al_k)$,  for $r\geq a_\eps'$, 
\[\F_{p_k,r}(F_k^\pm-C_{p_k}(\al^+_k))\geq\frac{\pi-\eps}{6}.\tag{$\star$}\] 
Since $\length(\al_k^+)\to\pi$, we have the uniform  mass bound 
\[(\|C_{p_k}(\al^+_k)\|+\|\D C_{p_k}(\al^+_k)\|)(B_r(p_k))\leq \pi r^2+2r.\]
Using \cite[Theorem~2.3]{KL_higher}, we may assume that $(\ga_k)_\#(C_{p_k}(\al^+_k))$ converges to a limit $C_\infty\in\bI_{2,loc}(X^+)$ with respect to local flat topology. By \cite[Proposition~2.2]{W_compact}, the support of $C_\infty$ is contained in $\tilde H$. 
Because $(\ga_k)_\#\D(C_{p_k}(\al^+_k))\to\D\tilde H$ weakly, we have $\D C_\infty=\D\tilde H$. But $C_\infty$
is a top-dimensional locally integral current in $\tilde H$, hence $C_\infty=\tilde H$ \cite[Theorem~7.2]{La_currents}.

Now suppose for contradiction, that $\tilde H=\tilde F^+$. Set $T_k=F_k^+-C_{p_k}(\al^+_k)\in\bZ_{2,loc}(X^+)$. Because $(\ga_k)_\#F_k^+\to\tilde F^+$ locally flat, we infer $(\ga_k)_\#(T_k)\to 0$ locally flat. Since $\D T_k=0$ for all $k\in\N$,
\cite[Proposition~3.2]{LaWe_loc} implies $\F_{p_\infty,r}((\ga_k)_\#(T_k))\to 0$ for every $r>0$.
Since $\ga_k(p_k)\to p_\infty$, this contradicts ($\star$). Hence $\tilde H$ cannot be contained in $\tilde F$ and the proof is complete.
\qed

\bibliographystyle{alpha}
\bibliography{rr}

\Addresses

\end{document}